\documentclass[a4paper,12pt]{amsart}
\usepackage{amssymb,amscd,amsmath,a4wide,graphicx,stmaryrd,fullpage,setspace,microtype,xr,ulem,textcomp,tikz,enumitem,accents}
\usepackage[pagebackref=true]{hyperref}
\usepackage[all]{xy}
\usepackage[T1]{fontenc}
\usepackage{lmodern}
\title{Left properness of flows}

\author[P. Gaucher]{Philippe Gaucher}

\address{Universit\'e de Paris, CNRS, IRIF, F-75006, Paris, France}

\urladdr{http://www.irif.fr/{\~{}}gaucher} 



\swapnumbers

\newcommand{\B}{\mathcal{B}}
\newcommand{\C}{\mathcal{C}}
\newcommand{\D}{\mathcal{D}}
\newcommand{\K}{\mathcal{K}}

\newcommand{\W}{\mathcal{W}}
\newcommand{\F}{\mathcal{F}}

\newcommand{\de}{\partial}
\newcommand{\p}{\times}

\renewcommand{\P}{\mathbb{P}}

\newcommand\circled[1]{%
\tikz[baseline=(C.base)]\node[draw,circle,inner sep=0pt](C) {#1};\!
}

\newtheorem*{thmN}{Theorem}

\newtheorem{thm}{Theorem}[section]
\newtheorem{prop}[thm]{Proposition}
\newtheorem{lem}[thm]{Lemma}
\newtheorem{exa}[thm]{Example}
\newtheorem{cor}[thm]{Corollary}

\newtheorem{defn}[thm]{Definition}
\newtheorem{nota}[thm]{Notation}
\newtheorem{rem}[thm]{Remark}
\newcommand{\bd}{\begin{defn}}
	\newcommand{\ed}{\end{defn}}
\newcommand{\bp}{\begin{prop}}
	\newcommand{\ep}{\end{prop}}
\newcommand{\bth}{\begin{thm}}
	\renewcommand{\eth}{\end{thm}}
\newcommand{\bpf}{\begin{proof}}
	\newcommand{\epf}{\end{proof}}
\newcommand{\bc}{\begin{cor}}
	\newcommand{\ec}{\end{cor}}

\newcommand{\fL}[1]{\ar@{->}[ll]_-{#1}}
\newcommand{\fR}[1]{\ar@{->}[rr]^-{#1}}
\newcommand{\fRr}[1]{\ar@{->}[rrr]^-{#1}}
\newcommand{\fD}[1]{\ar@{->}[dd]_-{#1}}
\newcommand{\fU}[1]{\ar@{->}[uu]^-{#1}}
\newcommand{\f}[2]{\ar@{->}[#1]|{#2}}
\newcommand{\ff}[2]{\ar@2{->}[#1]|{#2}}
\newcommand{\frr}[1]{\ar@{->}[rrrr]^-{#1}}

\newcommand{\fl}[1]{\ar@{->}[l]_-{#1}}
\newcommand{\fr}[1]{\ar@{->}[r]^-{#1}}
\newcommand{\fd}[1]{\ar@{->}[d]_-{#1}}
\newcommand{\fu}[1]{\ar@{->}[u]^-{#1}}

\renewcommand{\top}{{\mathbf{Top}}}
\newcommand{\Htop}{{\mathcal{H}\mathbf{Top}}}

\newcommand{\iso}{\cong}
\newcommand{\lp}{\left(}
\newcommand{\rp}{\right)}

\renewcommand{\leq}{\leqslant}
\renewcommand{\geq}{\geqslant}

\newcommand{\ptop}[1]{{\brm{{#1}dTop}}}

\def\cartesien{%
	\ar@{-}[]+R+<6pt,-2pt>;[]+RD+<6pt,-6pt>%
	\ar@{-}[]+D+<2pt,-6pt>;[]+RD+<6pt,-6pt>%
}
\def\cocartesien{%
	\ar@{-}[]+L+<-6pt,+2pt>;[]+LU+<-6pt,+6pt>%
	\ar@{-}[]+U+<-2pt,+6pt>;[]+LU+<-6pt,+6pt>%
}
\def\hocartesien{%
	\ar@{-}[]+R+<6pt,-2pt>;[]+RD+<6pt,-6pt>_{h}%
	\ar@{-}[]+D+<2pt,-6pt>;[]+RD+<6pt,-6pt>%
}
\def\hococartesien{%
	\ar@{-}[]+L+<-6pt,+2pt>;[]+LU+<-6pt,+6pt>_{h}%
	\ar@{-}[]+U+<-2pt,+6pt>;[]+LU+<-6pt,+6pt>%
}

\newcommand{\brm}[1]{{\rm{\mathbf{#1}}}}

\newcommand{\dtop}{{\brm{Flow}}}

\newcommand{\set}{{\brm{Set}}}

\newcommand{\ttop}{{\brm{TOP}}}
\newcommand{\mtop}{{\brm{MSpc}}}

\newcommand{\glob}{{\mathrm{Glob}}}

\newcommand{\globP}{{\mathrm{Glob}^{\mathcal{G}}}}

\DeclareMathOperator{\id}{Id}

\DeclareMathOperator{\Obj}{Obj}
\DeclareMathOperator{\Mor}{Mor}
\DeclareMathOperator{\pr}{pr}

\newcommand{\liminj}{\varinjlim}
\newcommand{\limproj}{\varprojlim}
\newcommand{\dcat}{{\mathbf{Cat}}}
\newcommand{\Dcat}{{\mathbf{CAT}}}

\newcommand{\rest}{\!\upharpoonright\!}

\renewcommand{\P}{\mathbb{P}}

\DeclareMathOperator{\Lan}{Lan}

\DeclareMathOperator{\CC}{CC}

\makeatletter
\def\varholim@#1#2{%
	\vtop{\m@th\ialign{##\cr
			\hfil$#1\operator@font holim$\hfil\cr
			\noalign{\nointerlineskip\kern1.5\ex@}#2\cr
			\noalign{\nointerlineskip\kern-\ex@}\cr}}%
}
\def\holimproj{%
	\mathop{\mathpalette\varholim@{\leftarrowfill@\textstyle}}\nmlimits@
}
\def\holiminj{%
	\mathop{\mathpalette\varholim@{\rightarrowfill@\textstyle}}\nmlimits@
}
\makeatother

\DeclareMathOperator{\cell}{{\brm{cell}}}
\DeclareMathOperator{\cof}{{\brm{cof}}}
\DeclareMathOperator{\inj}{{\brm{inj}}}
\newcommand{\ddownarrow}{{\downarrow}}

\DeclareMathOperator{\cocyl}{{Path}}
\newcommand{\adj}[5]{\xymatrix@C=#5em{{#1}\ar@/^0.8em/[r]^-{#2} \ar@{}[r]|-{\perp} & \ar@/^0.8em/[l]^-{#3} {#4}}}

\newcommand{\Di}{\mathbf{D}}
\newcommand{\Sp}{\mathbf{S}}
\newcommand{\De}{\Delta}

\begin{document}

\begin{abstract} 
Using Reedy techniques, this paper gives a correct proof of the left properness of the q-model structure of flows. It fixes the preceding proof which relies on an incorrect argument. The last section is devoted to fixing some arguments published in past papers coming from this incorrect argument. These Reedy techniques also enable us to study the interactions between the path space functor of flows with various notions of cofibrations. The proofs of this paper are written to work with many convenient categories of topological spaces like the ones of $k$-spaces and of weakly Hausdorff $k$-spaces and their locally presentable analogues, the $\Delta$-generated spaces and the $\De$-Hausdorff $\Delta$-generated spaces.
\end{abstract}

\maketitle

\tableofcontents

\section{Introduction}

\subsection*{Presentation}

The primary motivation for introducing the formalism of flows in \cite{model3} is the study of the branching areas and merging areas of execution paths in a concurrent system by using new homology theories. The idea of such homology theories dates back to \cite{Prehistoire-homology}. However, Goubault-Jensen's construction is not invariant by the refinement of observation. The more there are subdivisions in the description of the process, the more there are elements in the homology theories for the same branching or merging area. A \textit{flow} is a short name for a \textit{small topologically enriched semicategory} or \textit{small topologically enriched non-unital category}. The objects of a flow represent the \textit{states} of the concurrent process and the space of morphisms between two states is the topological space of \textit{execution paths}, the topology modelling concurrency. The reason for using $0$-dimensional identity morphisms instead of $1$-dimensional identity morphisms like in the formalism of small topologically enriched categories is to obtain functorial constructions for the \textit{branching and merging homology theories} (see \cite[Section~20]{model3}).  It is then proved in \cite[Corollary~11.3]{3eme} that the branching and merging homology theories of flows are invariant with respect to the refinement of observation \cite[Corollary~11.3]{3eme}. The latter paper therefore proves that the branching and merging homology theories of flows repair Goubault-Jensen's construction of \cite{Prehistoire-homology}. The main technical tool to study these branching and merging homology theories is the model structure introduced in \cite{model3}. It is called now the \textit{q-model structure of flows} after \cite{QHMmodel}. 

Precubical sets are a standard geometric model of concurrency which is widely used in the literature \cite{MR3445318}. Every precubical set can be viewed as a flow.  However, the naive realization functor which takes a $n$-cube to the flow associated with the poset of its vertices crushes the hollow $n$-cubes for $n\geq 3$ by \cite[Theorem~7.1]{ccsprecub}. Therefore, it does not yield a well-behaved realization functor. Roughly speaking, this geometric phenomenon is due to the fact that the permutations are generated by the transpositions. By replacing the commutativity relations in the partial monoid of execution paths by homotopies thanks to the q-model structure, it is possible to define a \textit{well-behaved realization functor} from precubical sets to flows (see \cite[Definition~7.2 and Theorem~7.6]{ccsprecub}). Once again, the homotopical techniques to study this realization functor use the q-model structure of flows. The study of this well-behaved realization functor is continued in \cite{realisation}, in which the key technical tool is also the q-model structure of flows.

The q-model structure of flows makes it also possible to define and to prove the main properties of the \textit{underlying homotopy type} of a flow \cite[Section~6]{4eme}. In particular, it is proved in \cite[Theorem~9.1]{4eme} that the underlying homotopy type of a flow is invariant with respect to the refinement of the observation, which is the expected behavior. Indeed, the underlying homotopy type of a flow is, morally speaking, the underlying space of states of a flow after removing the execution paths. It is defined only up to homotopy, not up to homeomorphism. Morally speaking, a flow is a directed space over a homotopy type indeed.

The notion of fundamental category has been proved to be a relevant object for static analysis of concurrent programs \cite{MR3445318}. It turns out that the weak equivalences of the q-model structure of flows preserve the \textit{fundamental category} of a flow. The latter is defined by taking the small category associated with the semicategory whose objects are the states of the flow and the set of morphisms between two states are the path-connected components of the space of morphisms between these two states \cite[Definition~2]{MR2057412}. This notion of fundamental category is easy to calculate because it is a left adjoint functor: it is therefore colimit-preserving. The right adjoint consists of taking a small category to the associated semicategory with discrete spaces of execution paths. This property is due to the fact that $\Delta$-generated spaces are homeomorphic to the disjoint sum of their path-connected components. Moreover the notion of a fundamental category of a flow interacts very well with the simplicial structure of the q-model structure of flows~\footnote{The model structure is simplicial by \cite[Theorem~3.3.15]{realisation}}, which makes it an interesting subject of study for future papers. The fundamental category of a $d$-space in the sense of \cite{mg}, defined in a similar way \cite[Definition~2]{MR2057412}, behaves slightly differently. The fundamental category of a flow remains finite for a finite precubical set without loops. It is not the case by using the formalism of Grandis' $d$-space. In that case, it gives rise to an enormous object containing uncountably many objects and morphisms. For example, the fundamental category of the flow associated with a $1$-dimensional cube (which models the directed segment) has two objects and one nonidentity morphism. On the contrary, the fundamental category of the $d$-space associated with a $1$-dimensional cube is the poset $([0,1],\leq)$. In the latter case, it is therefore necessary to introduce various notions of component categories to shrink the fundamental category \cite{MR2057412} \cite{MR2350213}. After an example due to Dubut of a finite precubical set without loops having an infinite component category in the previous sense \cite[page~162]{DubutPhD}, another way to reduce the size of the fundamental category is even proposed in \cite{MR4071310}. The short discussion of this paragraph illustrates the main difference between the geometric model of flows and other models of the literature, e.g. \cite{mg} \cite{MR2545830}. The former is a multipointed formalism, i.e. equipped with a distinguished set of objects, like the formalism of simplicial sets, the latter are not.  

All the previous examples show the theoretical importance of the q-model category of flows even if this model category does not have enough weak equivalences. Indeed, they identify less concurrent processes than we would like. It is the reason why it remains to understand the behavior of some Bousfield localizations, the left one with respect to the refinement of observation, the right one with respect to the $n$-cubes. The  q-model category of flows is right proper because all objects are fibrant. To study left Bousfield localizations in the framework of combinatorial model categories, left properness seems to be required even if a recent work enables us to get rid of this hypothesis \cite{LeftBousfield_withoutleftproperness}. The left properness of the q-model structure of flows is also used to prove the invariance with respect to the refinement of observation of the branching and merging homology theories in \cite[Theorem~11.2]{3eme}, and of the underlying homotopy type in \cite[Theorem~9.1]{4eme}. Left properness is difficult to prove in this kind of model category because pushouts, and more generally colimits of flows, can freely generate new execution paths in the colimit. The main result of this paper is a correct proof of the following theorem:  

\begin{thmN} (Theorem~\ref{proof-left-proper-flow} which corrects \cite[Theorem~7.4]{2eme})
	The q-model structure of flows is left proper. 
\end{thmN}

Another proof corrected in this paper is the one of the following theorem. Its statement is quite strange at first sight because, at least in the framework of $\Delta$-generated spaces ($\Delta$-Hausdorff or not), the path space functor is a right Quillen adjoint from the q-model structure of flows to the q-model structure of topological spaces: 

\begin{thmN} (Theorem~\ref{complement-flow} which corrects \cite[Proposition~A.2]{3eme})
	The path space functor $\P:\dtop\to \top$ preserves q-cofibrancy. 
\end{thmN}

In fact, the purpose of this paper is twofold. The first one is to fix the proof of \cite[Proposition~15.1]{model3} in Theorem~\ref{calculation-pathspace-flow}. It leads to a correct proof of Theorem~\ref{proof-left-proper-flow} and Theorem~\ref{complement-flow}. Section~\ref{flaw}, after explaining carefully the issue in the proof of  \cite[Proposition~15.1]{model3}, is mostly devoted to fixing the proof of some theorems of \cite{model2} which are used to prove that the homotopy categories of the q-model structures of multipointed $d$-spaces and flows are equivalent \cite[Theorem~7.5]{mdtop}. These problems have no influence on the theory of multipointed $d$-spaces and of flows as it has been developed so far. They only change the proofs of some intermediate results which remain valid anyway. 

The second one is to introduce some material required for the study of the homotopy theory of \textit{Moore flows}. This paper belongs to a series of papers (the order of publication is not the order of writing). It starts with \cite{leftdetflow} and \cite{QHMmodel} which revisit the q-model structure of flows. It continues with \cite{dgrtop} which establishes some theorems about the homotopy theory of enriched diagrams of topological spaces. Then the series continues with this paper. And it is finally concluded with the two papers \cite{Moore1} \cite{Moore2}. The purpose of the pair of papers \cite{Moore1} \cite{Moore2} is to upgrade the categorical equivalence between the homotopy categories  of the q-model structures of multipointed $d$-spaces and flows to a zigzag of Quillen equivalences. Moore flows are small semicategories enriched over a semimonoidal category of enriched presheaves of spaces over a specific reparametrization category. By taking the one-object category as a reparametrization category, we recover the notion of flow of this paper. It means that most of the theorems involving only the semicategorical nature of flows can be generalized to Moore flows. It is the case for example for the Reedy constructions of this paper, and for Theorem~\ref{calculation-pathspace-flow} which is not formulated in the most general way. 

I actually discovered the flaw in the proof of \cite[Proposition~15.1]{model3} and in the proof of its consequence \cite[Theorem~15.2]{model3} precisely by working on Moore flows. I wanted to prove that a q-cofibrant Moore flow has a projective q-cofibrant enriched presheaf of execution paths (it is a generalization of Theorem~\ref{complement-flow}). For this reason, I had to generalize the proof of \cite[Theorem~15.2]{model3}. I then realized that the proof of \cite[Proposition~15.1]{model3} was not correct. It led me from one thing to another to the Reedy constructions of this paper and to the statement of Theorem~\ref{calculation-pathspace-flow}.

\subsection*{Outline of the paper}
\begin{itemize}[leftmargin=*]
	\item Section~\ref{space} is a reminder about topological spaces. We want both to establish the results of this paper in the locally presentable setting of ($\De$-Hausdorff) $\De$-generated spaces (to prepare the subsequent papers) and to fix some past papers written in the framework of weakly Hausdorff $k$-spaces. It is the reason why we work in this paper in a framework containing all these situations as particular cases. We prove some important facts about relative-$T_1$ maps. The latter are a generalization of the notion of closed $T_1$-inclusion.
	\item Section~\ref{a-reedy-cat} introduces the Reedy category which is the keystone of the paper. It has both fibrant constants and cofibrant constants. Only the first property matters for this paper. It enables us to encode the calculation of all new execution paths created by a pushout along a generating q-cofibration. It is defined by generators and relations. It is proved that it is a poset.
	\item Section~\ref{calcul_final_structure-section} starts by a reminder about flows and its q-model structure. Then it is expounded in Theorem~\ref{calculation-pathspace-flow} the calculation of the path space of a flow which is obtained as a pushout along a map of flows of the form $\glob(\de Z) \subset \glob(Z)$. There is no hypothesis made on the continuous map $\de Z\to Z$ in this section. The case of a pushout along the generating cofibrations $C:\varnothing\subset \{0\}$ and $R:\{0,1\}\to \{0\}$ is not treated here because it is trivial. It will be just mentioned in the core of the proof of Theorem~\ref{proof-left-proper-flow}.
	\item Then Section~\ref{proof-leftproper-section} is entirely devoted to the proof of Theorem~\ref{proof-left-proper-flow}. The theory developed here is used, plus the fact that the homotopy colimits of a diagram of spaces in the q-model structure and in the h-model structure have the same weak homotopy type. This section is concluded by exploring in Theorem~\ref{complement-flow} the interactions between the path space functor of flows and the classes of cofibrations of the model structures we have worked with. In particular, we find a proof that the path space functor from flows to topological spaces preserves q-cofibrancy. These interactions are surprising because the path space functor of flows is a right Quillen adjoint in the locally presentable case. Note that most of the results of Section~\ref{proof-leftproper-section} are new.
	\item The precise description of the flaw in the proof of \cite[Proposition~15.1]{model3} is postponed to Section~\ref{flaw}. Then we explain why \cite[Theorem~15.2]{model3} is true anyway despite the  incorrect argument. As for the group of papers \cite{3eme} \cite{4eme} \cite{2eme} \cite{mdtop}, it is explained how not to use the same (probably~\footnote{I cannot prove that the lemma is  wrong; however, I am sure that the argument leading to it is  wrong.})  wrong lemma coming from the flaw. Finally it is explained why \cite[Theorem~V.3.4]{model2} is still true (after removing an assertion which is useless and that it is not known whether it is true) in Theorem~\ref{fix} and why \cite[Theorem~III.5.2]{model2} is still true by supplying in Theorem~\ref{existence-section} an updated proof using the tools developed in this paper.
	\item Appendix~\ref{alldg} expounds some very basic properties about the category of all diagrams over all small categories valued in a bicomplete category. 
	\item Appendix~\ref{HTop} introduces a notion of separation on $\De$-generated spaces. This new setting is another convenient category of topological spaces for doing algebraic topology. It avoids dealing with pointless point set topology problems involving the indiscrete topology and the Sierpinski topology while preserving the local presentability of the underlying category of spaces. Indeed, all spaces are $T_1$ in this category. This appendix does not pretend to be exhaustive. It proves what is needed for the paper and nothing more. 
\end{itemize}

\subsection*{Notations and terminology}

We refer to \cite{TheBook} for locally presentable categories, to \cite{MR2506258} for combinatorial model categories.  We refer to \cite{MR99h:55031} and to \cite{ref_model2} for more general model categories. 

\begin{itemize}[leftmargin=*]
	\item $X:=Y$ means that $Y$ is the definition of $X$.
	\item All categories are locally small (except the category $\Dcat$ of all locally small categories).
	\item $\dcat$ is the category of small categories and functors between them.
	\item $\set$ is the category of sets.
	\item $\mathcal{T\!O\!P}$ is the category of general topological spaces with the continuous applications.
	\item A final quotient of $\mathcal{T\!O\!P}$ is a surjective continuous map $q:X\to Y$ such that $Y$ is equipped with the final topology.
	\item A map of $\mathcal{T\!O\!P}$ is always supposed to be continuous; otherwise the terminology \textit{set map} is used.
	\item The paper uses the French convention: compact implies Hausdorff. A topological space satisfying the finite open covering property is called quasi-compact.
	\item $\B$ denotes either the singleton $\mathbf{\De}=\{[0,1]\}$ where $[0,1]$ is the compact segment or the proper class $\mathbf{K}$ of all nonempty compact spaces. 
	\item $\top_\B$ is the final closure of $\B$ in $\mathcal{T\!O\!P}$. 
	\item $\Htop_\B$ is the full subcategory of $\top_\B$ of $\B$-Hausdorff spaces. 
	\item The inclusion functor $\top_\B\subset \mathcal{T\!O\!P}$ has a right adjoint: the $\B$-kelleyfication functor $k_\B:\mathcal{T\!O\!P}\to \top_\B$. 
	\item The inclusion functor $\Htop_\B\subset \top_\B$ has a left adjoint: the $\B$-Hausdorffization functor $w_\B:\top_\B\to \Htop_\B$.
	\item The category $\top$ is one of the categories $\top_\B$ or $\Htop_\B$ with $\B$ equal to $\mathbf{\De}$ or $\mathbf{K}$.
	\item $\K(X,Y)$ is the set of maps in a category $\K$ from $X$ to $Y$. $|\K|$ is the class of objects.
	\item A transfinite tower (of length $\lambda$) of $\K$ consists of a limit ordinal $\lambda$ and a colimit-preserving functor $D$ from $\lambda$ to $\K$; it means that for every limit ordinal $\mu\leq \lambda$, the canonical map $\liminj_{\nu<\mu} D_\nu\to D_\mu$ is an isomorphism.
	\item $A \sqcup B$ is the binary coproduct, $A \p B$ is the binary product. 
	\item $\limproj$ is the limit, $\liminj$ is the colimit.
	\item $\varnothing$ is the initial object.
	\item $\mathbf{1}$ is the terminal object.
	\item $\id_X$ is the identity of $X$.
	\item $g.f$ is the composite of two maps $f:A\to B$ and $g:B\to C$; the composite of two functors is denoted in the same way. 
	\item $\K^I$ is the category of functors and natural transformations from a small category $I$ to $\K$.
	\item If $f:{I}\to {J}$ is a functor between small categories and if $F:{I}\to\K$ is a functor, then $\Lan_f F$ denotes the left Kan extension of $F$ along $f$.
	\item $F\Rightarrow G$ denotes a natural transformation from a functor $F$ to a functor $G$.
	\item The composite of two natural transformations $\mu:F\Rightarrow G$ and $\nu:G\Rightarrow H$ is denoted by $\nu \odot \mu$ to make the distinction with the composition of maps.
	\item $f\boxslash g$ means that $f$ satisfies the \textit{left lifting property} (LLP) with respect to $g$, or equivalently that $g$ satisfies the \textit{right lifting property} (RLP) with respect to $f$.
	\item $\inj(\C) = \{g \in \K, \forall f \in \C, f\boxslash g\}$.
	\item $\cof(\C)=\{f\mid \forall g\in \inj(\C), f\boxslash g\}$.
	\item $\cell(\C)$ is the class of transfinite compositions of pushouts of elements of $\C$.
	\item A \textit{cellular} object $X$ of a cofibrantly generated model category is an object such that the canonical map $\varnothing\to X$ belongs to $\cell(I)$ where $I$ is the set of generating cofibrations.
	\item A \textit{model structure} $(\C,\W,\F)$ means that the class of cofibrations is $\C$, that the class of weak equivalences is $\W$ and that the class of fibrations is $\F$ in this order. A \textit{model category} is a category equipped with a \textit{model structure}.
	\item If $\D:I\to\K$ is a diagram over a Reedy category $(I,I_+,I_-)$, the latching category at $i\in I$ is denoted by $\de(I_+\ddownarrow i)$,  the latching object at $i\in I$ by $L_i\D := \liminj_{\de(I_+\ddownarrow i)} \D$, the matching category at $i\in I$ by $\de(i\ddownarrow I_-)$ and the matching object at $i\in I$ by $M_i\D = \limproj_{\de(i\ddownarrow I_-)} \D$.
	\item $f\square g$ is the pushout product of two maps $f$ and $g$. 
	\item $\pi_n(X)$ means the $n$-th homotopy group of $X$ for some base point.
	\item A cocone from a diagram $\D$ to an object $X$ is denoted by $\D\stackrel{\bullet}\to X$. 
	\item The $n$-dimensional disk for $n\geq 1$ is denoted by $\Di^n = \{(x_1,\dots,x_n)\in \mathbb{R}^n\mid x_1^2 + \dots + x_n^2 \leq 1\}$. By convention, let $\Di^0=\{0\}$. 
	\item The $(n-1)$-dimensional sphere for $n\geq 1$ is denoted by $\Sp^n = \{(x_1,\dots,x_n)\in \mathbb{R}^n\mid x_1^2 + \dots + x_n^2 = 1\}$. By convention, let $\Sp^{-1}=\varnothing$.
	\item The $n$-dimensional simplex for $n\geq 0$ is denoted by $\Delta^n = \{(x_0,\dots,x_n)\in [0,1]^{n+1}\mid x_0 + \dots + x_n = 1\}$.
	\item All h-cofibrations are by convention \textit{strong} h-cofibrations in the sense of \cite{MR1967263} \cite{Barthel-Riel}. The terminology of $\overline{h}$-cofibration and  $\overline{h}$-fibration is not used.
	\item An inclusion $i:A\to B$ is a one-to-one map such that $A$ is homeomorphic to $i(A)$ equipped with the relative topology.
	\item A space is discrete if it is equipped with the discrete topology. A space is totally disconnected if its connected components are its points.
\end{itemize}

\subsection*{Acknowledgment}

I thank the anonymous referee for the report and for the helpful remarks to improve the introduction.

\section{Topological spaces}
\label{space}

A \textit{$\B$-generated space} is a topological space which belongs to the final closure of $\B$ in $\mathcal{T\!O\!P}$. A general topological space $X$ is \textit{$\B$-Hausdorff} if for every continuous map $f:K\to X$ with $K\in \B$, the set $f(K)$ is a closed subset of $X$. The case $\B=\mathbf{K}$ is well-known. The case $\B=\mathbf{\De}$ is treated in Appendix~\ref{HTop}. The reason for working at this level of generality is that Section~\ref{flaw} is devoted to fixing some past proofs written in the category of weakly Hausdorff $k$-spaces. We summarize first some basic properties of $\top$~\footnote{which is one of the categories $\top_\B$ or $\Htop_\B$ with $\B$ equal to $\mathbf{\De}$ or $\mathbf{K}$.} needed for this work for the convenience of the reader: 
\begin{itemize}
	\item The $\B$-kelleyfication functor $k_\B:\mathcal{T\!O\!P} \to \top_\B$ does not change the underlying set.
	\item Let $A\subset B$ be a subset of a space $B$ of $|\top|$. Then $A$ equipped with the $\B$-kelleyfication of the relative topology belongs to $|\top|$. Note that for $\B=\mathbf{\De}$, a closed subset of a $\B$-generated space equipped with the relative topology is not necessarily $\B$-generated: e.g. the Cantor set $\mathbb{K}\subset [0,1]$ is not $\Delta$-generated; its $\Delta$-kelleyfication is the Cantor set equipped with the discrete topology $\mathbb{K}^\delta$. It is always the case if $\B=\mathbf{K}$. An open subset of a $\Delta$-generated space equipped with the relative topology is always $\Delta$-generated. This comes from the fact that any open subset of $[0,1]$ is $\De$-generated (see also \cite[Proposition~1.18]{delta}).
	\item $\top$ is cartesian closed. The internal hom $\ttop(X,Y)$ is given by taking the $\B$-kelleyfication of the compact-open topology on the set $\mathcal{T\!O\!P}(X,Y)$ of all continuous maps from $X$ to $Y$. 
	\item The colimit in $\top$ is given by the final topology in the following situations: 
	\begin{itemize}
		\item A transfinite compositions of one-to-one maps.
		\item A pushout along a closed inclusion.
		\item A quotient by a closed subset or by an equivalence relation having a closed graph.
	\end{itemize}
	In these cases, the underlying set of the colimit is therefore the colimit of the underlying sets. In particular, the CW-complexes, and more generally all cellular spaces are equipped with the final topology. Note that cellular spaces are even Hausdorff (and paracompact, normal, etc...). 
	\item The category $\top$ admits a q-model structure, a h-model structure and a m-model structure. All q-cofibrations are m-cofibrations and all m-cofibrations are h-cofibrations. 
\end{itemize}

\begin{rem}
	All h-cofibrations are by convention \textit{strong} h-cofibrations in the sense of \cite{MR1967263} \cite{Barthel-Riel}. The terminology of $\overline{h}$-cofibration and  $\overline{h}$-fibration is not used. It means that a h-cofibration is a closed inclusion satisfying the LLP with respect to all maps of the form $\ttop([0,1],Y)\to Y$. 
\end{rem}

Both $\Htop_{\mathbf{\De}}$ and $\top_{\mathbf{\De}}$ are locally presentable and every $\De$-generated space is homeomorphic to the disjoint sum of its nonempty path-connected components which are also its nonempty connected components. The latter hypothesis is used only in Theorem~\ref{rightadjointpath-flow}. The local presentability is not used in the core of the paper.

Both $\Htop_\mathbf{K}$ and $\top_\mathbf{K}$ have a h-model structure by \cite[Corollary 5.23]{Barthel-Riel}: they satisfy the monomorphism hypothesis by \cite[Example~5.18]{Barthel-Riel} and they are topologically bicomplete because they are cartesian closed.

The following notion is a weakening of the notion of closed $T_1$-inclusion introduced by Dugger and Isaksen. It enables us to work with or without the separation condition (like in \cite{mdtop}).

\bd \cite[p 686]{hocolimfacile} A one-to-one continuous map $i:A\rightarrow X$ is {\rm relative-$T_1$} if for any open subset $U$ of $A$ and any point $z\in X\backslash i(U)$, there is an open set $W$ of $X$ with $i(U)\subset W$ and $z\notin W$. \ed

We have: 

\bp
Let $i:A\to X$ be a closed $T_1$-inclusion. Then $i$ is relative-$T_1$. 
\ep

\bpf
Let $U$ be an open subset of $A$. Let $z\in X\backslash i(U)$. Assume that $z\in i(A)$, and write $i(U)=i(A)\cap W$ for some open subset $W$ of $X$. Then $\{z\}\cap W=(\{z\} \cap i(A)) \cap W = \{z\} \cap (i(A) \cap W)=\{z\} \cap i(U)=\varnothing$. We deduce that $z\notin W$ and $i(U)=i(A)\cap W \subset W$. Now assume that $z\notin i(A)$. Then let $W=X\backslash \{z\}$. Since $z$ is a closed point by hypothesis, $W$ is an open subset of $X$. It does not contain $z$ and $i(U)\subset i(A) \subset W$.  
\epf

The following proposition gives an example of a relative-$T_1$ inclusion which is neither closed nor $T_1$.

\bp
Let $X$ be a $\De$-generated space which is not $T_1$ (for example, an indiscrete space). Let $A$ be an open subset of $X$ equipped with the relative topology. Then $A$ is a $\De$-generated space and the inclusion $A\subset X$ is a relative-$T_1$ inclusion.
\ep

\bpf
Let $U$ be an open of $A$ and $z\in X\backslash U$. Then $W=U\cap A=U$ is an open subset of $X$ and $z\notin W$.
\epf

The following proposition is a replacement of the usual one: 

\bp \label{factor-tower} Every final quotient of a space of $\B$ is finite relative to relative-$T_1$ inclusions. 
\ep

\bpf
Since the colimit in $\top$ of a tower of one-to-one maps is always equipped with the final topology, the proposition is a consequence of \cite[Lemma~A.3]{hocolimfacile}. 
\epf

There is the key fact:

\bp \label{h-cof-relative-T1} All h-cofibrations of $\top$ are relative-$T_1$. \ep

\bpf
Let $i:A \to X$ be a h-cofibration of $\top$. Then it is a closed inclusion. Therefore we can suppose that $A\subset X$ with $A$ equipped with the relative topology. The rest of the proof is a modification of the one of \cite[Proposition~1(b)]{vstrom3}. Consider \[E=A\p [0,1]  \cup X \p ]0,1] \subset X \p [0,1].\] Consider a commutative diagram of spaces of the form
\[
\xymatrix@C=4em@R=4em
{
	Y\p \{0\} \fd{\subset} \fr{g} & E \ar@{->}[d]^-{p} \fr{\pr_{[0,1]}}  & [0,1]\\
	Y\p [0,1] \ar@{-->}[ru]^-{\overline{G}}  \fr{G} & X .
}
\] 
where $p:E\to X$ is the projection map and where $\pr_{[0,1]}$ is the projection over $[0,1]$. Then $\overline{G}(y,t)=(G(y,t),t+(1-t)\pr_{[0,1]}g(y))$ is a lift. It means that the map $p:E\to X$ is a h-fibration. Consider the commutative diagram of spaces, with $k(a)=(a,0)$:
\[
\xymatrix@C=4em@R=4em
{
	A \ar@{->}[d]_-{i} \fr{k} & E \fr{\pr_{[0,1]}} \ar@{->}[d]^-{p} & [0,1]\\
	X \ar@{=}[r] \ar@{-->}[ru]^-{f} &  X .
}
\]
Using the homotopy $H:((x,t),u)\mapsto (x,u+(1-u)t)$, we see that $p$ is a homotopy equivalence. Therefore the lift $f$ exists because, by hypothesis, the map $i:A\to X$ is a h-cofibration. Let $U$ be an open subset of $A$ and let $z\in X\backslash U$. There are two mutually exclusive cases: 
\begin{enumerate}
	\item $z\in A$. We have $U=A \cap W$ for some open subset $W$ of $X$. Then $\{z\}\cap W = (\{z\}\cap A) \cap W = \{z\}\cap (A \cap W) = \{z\}\cap U = \varnothing$. Thus $z\notin W$ and $U=A \cap W\subset W$.
	\item $z\notin A$. Then $\pr_{[0,1]}(f(z)) \in ]0,1]$ because $A = (f^{-1}.\pr_{[0,1]}^{-1})(0)$. Consider the open subset $W=f^{-1}(\pr_{[0,1]}^{-1})([0,\pr_{[0,1]}(f(z))/2[)$ of $X$. Then by construction, $U\subset A = (f^{-1}.\pr_{[0,1]}^{-1})(0)\subset W$ and $z\notin W$.  
\end{enumerate}
We have proved that $i:A \to X$ is relative-$T_1$.
\epf

We obtain the important consequence: 

\begin{cor} Every final quotient in $\top$ of a space of $\B$ is finite relative to h-cofibrations. \end{cor}

We conclude this section with another theorem which plays a central role in this work:

\bth (Dugger-Isaksen) \label{hocolimfacile1} Let $D:I\to \top$ be a small diagram. Then the homotopy colimits of $D$ as computed in the q-model structure of  $\top$ and in the h-model structure of  $\top$ have the same weak homotopy type. \eth

\bpf[Sketch of proof] We explain the difference with the proof of \cite[Theorem~A.7]{hocolimfacile} written in the category of general topological spaces. 
\begin{itemize}
	\item \cite[Lemma~A.1]{hocolimfacile} still holds in  $\top$. Indeed, if $A$ and $B$ are $\B$-Hausdorff, then $A\p \Di^n$, $A\p \mathbf{S}^{n-1}$, $B\p \Di^n$ and $B\p \mathbf{S}^{n-1}$ with the product taken in  $\top_{\B}$ are $\B$-Hausdorff since  $\Htop_{\B}$ is a reflective subcategory of  $\top_{\B}$. Moreover the maps $A\p \mathbf{S}^{n-1} \subset A\p \Di^n$ and $B\p \mathbf{S}^{n-1} \subset B\p \Di^n$ are closed inclusions. Thus $X_A$ and $Y_B$ equipped with the final topology are $\B$-Hausdorff. It means that the underlying set of all spaces involved in the proof of \cite[Lemma~A.1]{hocolimfacile} are the same.
	\item \cite[Lemma~A.2]{hocolimfacile} is still valid in $\top$ without change.
	\item \cite[Lemma~A.3]{hocolimfacile} holds for any final quotient of $[0,1]$ (see Proposition~\ref{factor-tower}), so in particular for the $n$-spheres which is how it is used in \cite{hocolimfacile}.
	\item Then \cite[Proposition~A.5]{hocolimfacile} and \cite[Corollay~A.6]{hocolimfacile} follow.
	\item \cite[Theorem~A.7]{hocolimfacile} is precisely the statement of the theorem. 
\end{itemize}
\epf

\section{A Reedy category}
\label{a-reedy-cat}

Let $S$ be a nonempty set. Let $\mathcal{P}^{u,v}(S)$ be the small category defined by generators and relations as follows: 
\begin{itemize}[leftmargin=*]
	\item $u,v\in S$ ($u$ and $v$ may be equal).
	\item The objects are the tuples of the form 
	\[\underline{m}=((u_0,\epsilon_1,u_1),(u_1,\epsilon_2,u_2),\dots ,(u_{n-1},\epsilon_n,u_n))\]
	with $n\geq 1$, $u_0,\dots,u_n \in S$, $\epsilon_1,\dots,\epsilon_n \in \{0,1\}$ and \[\forall i\hbox{ such that } 1\leq i\leq n, \epsilon_i = 1\Rightarrow (u_{i-1},u_i)=(u,v).\] The integer $n$ is the \textit{length} of the tuple. The integer $\sum_i \epsilon_i$ is the \textit{height} of the tuple. 
	\item There is an arrow \[c_{n+1}:(\underline{m},(x,0,y),(y,0,z),\underline{n}) \to (\underline{m},(x,0,z),\underline{n})\]
	for every tuple $\underline{m}=((u_0,\epsilon_1,u_1),(u_1,\epsilon_2,u_2),\dots ,(u_{n-1},\epsilon_n,u_n))$ with $n\geq 1$ and every tuple $\underline{n}=((u'_0,\epsilon'_1,u'_1),(u'_1,\epsilon'_2,u'_2),\dots ,(u'_{n'-1},\epsilon'_{n'},u'_{n'}))$ with $n'\geq 1$. It is called a \textit{composition map}. 
	\item There is an arrow \[I_{n+1}:(\underline{m},(u,0,v),\underline{n}) \to (\underline{m},(u,1,v),\underline{n})\] for every tuple $\underline{m}=((u_0,\epsilon_1,u_1),(u_1,\epsilon_2,u_2),\dots ,(u_{n-1},\epsilon_n,u_n))$ with $n\geq 1$ and every tuple $\underline{n}=((u'_0,\epsilon'_1,u'_1),(u'_1,\epsilon'_2,u'_2),\dots ,(u'_{n'-1},\epsilon'_{n'},u'_{n'}))$ with $n'\geq 1$.
	It is called an \textit{inclusion map}. 
	\item There are the relations (group A) $c_i.c_j = c_{j-1}.c_i$ if $i<j$ (which means since $c_i$ and $c_j$ may correspond to several maps that if $c_i$ and $c_j$ are composable, then there exist $c_{j-1}$ and $c_i$ composable satisfying the equality). 
	\item There are the relations (group B) $I_i.I_j = I_j.I_i$ if $i\neq j$. By definition of these maps, $I_i$ is never composable with itself. 
	\item There are the relations (group C) \[c_i.I_j = \begin{cases}
	I_{j-1}.c_i&\hbox{if } j\geq i+2\\
	I_j.c_i&\hbox{if } j\leq i-1.
	\end{cases}\]
	By definition of these maps, $c_i$ and $I_i$ are never composable as well as $c_i$ and $I_{i+1}$. 
\end{itemize}

\bd Denote by $\mathcal{P}^{u,v}(S)_+$ the subcategory of $\mathcal{P}^{u,v}(S)$ generated by all objects of $\mathcal{P}^{u,v}(S)$ and by the inclusion maps. Denote by $\mathcal{P}^{u,v}(S)_-$ the subcategory of $\mathcal{P}^{u,v}(S)$ generated by all objects of $\mathcal{P}^{u,v}(S)$ and by the composition maps.
\ed

\bp \label{posetmoins} The category $\mathcal{P}^{u,v}(S)_-$ is a poset.
\ep

\bpf
Suppose that $\mathcal{P}^{u,v}(S)_-(\underline{m},\underline{n})$ is nonempty. Then, by definition of the composition maps, there are the equalities $\underline{m}=((u_0,\epsilon_1,u_1),(u_1,\epsilon_2,u_2),\dots ,(u_{n-1},\epsilon_n,u_n))$ and $\underline{n}=((u'_0,\epsilon'_1,u'_1),(u'_1,\epsilon'_2,u'_2),\dots ,(u'_{n'-1},\epsilon'_{n'},u'_{n'}))$ with $1\leq n'\leq n$ and the tuple $(u'_1,\dots,u'_{n'})$ is obtained from the tuple $(u_1,\dots,u_{n})$ by removing $u_{i_1},\dots,u_{i_p}$ with $i_1<\dots < i_p$ with $p\geq 0$. Using the presimplicial relations of group A, we see that the unique map of $\mathcal{P}^{u,v}(S)_-(\underline{m},\underline{n})$ is the composite map $c_{i_1}\dots c_{i_p}$. Because composition maps decrease the length, there is no pair of distinct isomorphic objects and the small category $\mathcal{P}^{u,v}(S)_-$ is a poset.
\epf

\bd An object $\underline{n}$ of the small category $\mathcal{P}^{u,v}(S)$ is {\rm simplifiable} if the matching category $\de(\underline{n}\ddownarrow \mathcal{P}^{u,v}(S)_-)$ is nonempty. 
\ed

\bp \label{final-matching}
Let $\underline{n}$ be an object of $\mathcal{P}^{u,v}(S)$. Then either $\underline{n}$ is not simplifiable (in this case, let $\mathbb{S}(\underline{n}):=\underline{n}$) or the matching category $\de(\underline{n}\ddownarrow \mathcal{P}^{u,v}(S)_-)$ has a terminal object denoted by $\mathbb{S}(\underline{n})$ and the latter is not simplifiable. 
\ep

\bpf
Let $\underline{n} = ((u_0,\epsilon_1,u_1),(u_1,\epsilon_2,u_2),\dots ,(u_{n-1},\epsilon_n,u_n))$. A descending chain of maps $\underline{n} \to \bullet \to \bullet \to\dots$ of $\mathcal{P}^{u,v}(S)_-$ will stop eventually since the length decreases along the chain. Let $\mathbb{S}(\underline{n})$ be a target of a maximal descending chain of maps of $\mathcal{P}^{u,v}(S)_-$. Then, by definition of the composition maps, we have necessarily 
\[
\mathbb{S}(\underline{n}) = ((u_{i_0},\epsilon'_1,u_{i_1}),(u_{i_1},\epsilon'_2,u_{i_2}),\dots ,(u_{i_{n'-1}},\epsilon'_{n'},u_{i_{n'}}))
\]
with $0=i_0 <i_1<\dots <i_{n'}=n$ with
\[\epsilon'_j = \begin{cases}
0 &\hbox{if } i_j - i_{j-1} > 1\\
\epsilon_{i_j} &\hbox{if } i_j - i_{j-1} = 1,
\end{cases}\]
with never two consecutive zeros in the sequence $\epsilon'_1,\dots, \epsilon'_{n'}$. Let \[\mathbb{S}'(\underline{n}) = ((u_{i'_0},\epsilon''_1,u_{i'_1}),(u_{i'_1},\epsilon''_2,u_{i'_2}),\dots ,(u_{i'_{n''-1}},\epsilon''_{n''},u_{i_{n''}}))\] be another target of a maximal descending chain of maps of $\mathcal{P}^{u,v}(S)_-$. Then $i_0=i'_0=0$. We necessarily have $\epsilon'_1=\epsilon''_1$. If $\epsilon'_1=\epsilon''_1=1$, then $i_1=i'_1=1$. If $\epsilon'_1=\epsilon''_1=0$ and e.g. $i_1<i'_1$, then $\mathbb{S}(\underline{n})$ is not a target of a maximal descending chain: contradiction. We deduce that $i_1=i'_1$. Proceeding by induction, we deduce that $\mathbb{S}(\underline{n})=\mathbb{S}'(\underline{n})$ is unique. Using Proposition~\ref{posetmoins}, the proof is complete. 
\epf

\bp \label{posetplus} The category $\mathcal{P}^{u,v}(S)_+$ is a poset. 
\ep

\bpf 
Suppose that $\mathcal{P}^{u,v}(S)_+(\underline{m},\underline{n})$ is nonempty. Then, by definition of the inclusion maps, there are the equalities $\underline{m}=((u_0,\epsilon_1,u_1),(u_1,\epsilon_2,u_2),\dots ,(u_{n-1},\epsilon_n,u_n))$ and $\underline{n}=((u_0,\epsilon'_1,u_1),(u_1,\epsilon'_2,u_2),\dots ,(u_{n-1},\epsilon'_{n'},u_{n}))$ with $n\geq 1$ and for all $i=1,\dots,n$, there is the inequality $\epsilon_i\leq \epsilon'_i$. Using the relations of group B, we see that the unique map of $\mathcal{P}^{u,v}(S)_+(\underline{m},\underline{n})$ is the composition $\prod_{\{i\mid \epsilon_i<\epsilon'_i\}} I_i$. Because inclusion maps increase the height, there is no pair of distinct isomorphic objects and the small category $\mathcal{P}^{u,v}(S)_+$ is a poset.
\epf

\bp \label{initial-latching}
Let $\underline{n}$ be an object of $\mathcal{P}^{u,v}(S)$. Then either $\de(\mathcal{P}^{u,v}(S)_+\ddownarrow\underline{n})$ is empty (in this case, let $\mathbb{I}(\underline{n}):=\underline{n}$) or it has an initial object denoted by $\mathbb{I}(\underline{n})$. 
\ep

\bpf
Let $\underline{n} = ((u_0,\epsilon_1,u_1),(u_1,\epsilon_2,u_2),\dots ,(u_{n-1},\epsilon_n,u_n))$. Then we have necessarily 
\[
\mathbb{I}(\underline{n}) = ((u_0,0,u_1),(u_1,0,u_2),\dots ,(u_{n-1},0,u_n)).
\]
The proposition is then a consequence of Proposition~\ref{posetplus}.
\epf

\bp \label{reedy}
The pair $(\mathcal{P}^{u,v}(S)_+,\mathcal{P}^{u,v}(S)_-)$ endows the small category $\mathcal{P}^{u,v}(S)$ with a structure of Reedy category with the $\mathbb{N}$-valued degree map defined by \[d((u_0,\epsilon_1,u_1),(u_1,\epsilon_2,u_2),\dots ,(u_{n-1},\epsilon_n,u_n)) = n + \sum_i \epsilon_i.\]
Moreover, in the canonical decomposition $f=f_+.f_-$ with $f_+ \in \Mor(\mathcal{P}^{u,v}(S)_+)$ and $f_-\in \Mor(\mathcal{P}^{u,v}(S)_-)$, the source of $f_+$, which is the target of $f_-$, is uniquely determined by the source and the target of $f$.
\ep

The minimal value of the degree map is $1$ and it is reached for the objects $((u_0,0,u_1))$ for $(u_0,u_1)$ running over $S\p S$. 

\bpf
The composition maps decrease the degree by one, the inclusion maps increase the degree by one. So every map of $\mathcal{P}^{u,v}(S)_+$ increases the degree and every map of $\mathcal{P}^{u,v}(S)_-$ decreases the degree. Let 
\begin{multline*}
f:((u_0,\epsilon_1,u_1),(u_1,\epsilon_2,u_2),\dots ,(u_{n-1},\epsilon_n,u_n))\\\to ((u'_0,\epsilon'_1,u'_1),(u'_1,\epsilon'_2,u'_2),\dots ,(u'_{n-1},\epsilon'_{n'},u'_{n'}))
\end{multline*}
be a map of $\mathcal{P}^{u,v}(S)$. By definition of the small category $\mathcal{P}^{u,v}(S)$, $f$ is a composite of composition maps and of inclusion maps. Using the relations of group C, we obtain a factorization  $f=f_+.f_-$ with $f_+\in \Mor(\mathcal{P}^{u,v}(S)_+)$ and $f_-\in \Mor(\mathcal{P}^{u,v}(S)_-)$. By definition of the inclusion maps, the source of $f_+$, which is the target of $f_-$, is of the form 
\[
((u'_0,\epsilon''_1,u'_1),(u'_1,\epsilon''_2,u'_2),\dots ,(u'_{n-1},\epsilon''_{n'},u'_{n'}))
\]
with $\epsilon''_j\leq  \epsilon'_j$ for $1\leq j\leq n'$. And by definition of the composition maps, there is the equality $(u'_0,u'_1,\dots,u'_{n'}) = (u_{i_0},u_{i_1},\dots,u_{i_{n'}})$ where $0= i_0 < i_1<\dots < i_{n'}= n$ and with
\[\epsilon''_j = \begin{cases}
 0 &\hbox{if } i_j - i_{j-1} > 1\\
 \epsilon_{i_j} &\hbox{if } i_j - i_{j-1} = 1.
\end{cases}\]
In other terms, there is only one possibility for the source of $f_+$ which is the target of $f_-$. The proof is complete thanks to Proposition~\ref{posetmoins} and Proposition~\ref{posetplus}.
\epf

\begin{cor} \label{poset}
The small category $\mathcal{P}^{u,v}(S)$ is a poset.
\end{cor}

We could directly define $\mathcal{P}^{u,v}(S)$ as a poset. The interest of having a presentation by generators and relations is that the proof of Proposition~\ref{well-defined-diag} becomes trivial. We will use this Reedy category as follows: 

\bth \label{ok} Let $\K$ be a model category. Let $S$ be a nonempty set. Let $u,v\in S$. Let $\Dcat(\mathcal{P}^{u,v}(S),\K)$ be the category of functors and natural transformations from $\mathcal{P}^{u,v}(S)$ to $\K$. Then there exists a unique model structure on $\Dcat(\mathcal{P}^{u,v}(S),\K)$ such that the weak equivalences are the pointwise weak equivalences and such that a map of diagrams $f:\D\to \mathcal{E}$ is a cofibration (called a Reedy cofibration) if for all objects $\underline{n}$ of $\mathcal{P}^{u,v}(S)$, the canonical map $L_{\underline{n}}\mathcal{E} \sqcup_{L_{\underline{n}} \D} \D(\underline{n})\to \mathcal{E}(\underline{n})$ is a cofibration of $\K$. Moreover the colimit functor $\liminj : \Dcat(\mathcal{P}^{u,v}(S),\K) \to \K$ is a left Quillen adjoint. \eth

\bpf A model structure is characterized by its class of weak equivalences and its class of cofibrations. Hence the uniqueness. The existence is explained e.g. in \cite[Theorem~15.3.4]{ref_model2}. The matching category of an object is either empty or connected by Proposition~\ref{final-matching}. The last assertion is then the consequence of \cite[Proposition~15.10.2]{ref_model2} and \cite[Theorem~15.10.8]{ref_model2}.
\epf

Note that the limit functor $\limproj : \Dcat(\mathcal{P}^{u,v}(S),\K) \to \K$ is a right Quillen adjoint by Proposition~\ref{initial-latching}, \cite[Proposition~15.10.2]{ref_model2} and \cite[Theorem~15.10.8]{ref_model2}.

\section{Path space of a pushout of flows along a {q}-cofibration}
\label{calcul_final_structure-section}

\bd \cite{model3} A {\rm flow} $X$ consists of a topological space $\P X$ of execution paths, a discrete space $X^0$ of states, two continuous maps $s$ and $t$ from $\P X$ to $X^0$ called the source and target map respectively, and a continuous and associative map \[*:\{(x,y)\in \P X\p \P X; t(x)=s(y)\}\longrightarrow \P X\] such that $s(x*y)=s(x)$ and $t(x*y)=t(y)$.  A morphism of flows $f:X\longrightarrow Y$ consists of a set map $f^0:X^0\longrightarrow Y^0$ together with a continuous map $\P f:\P X\longrightarrow \P Y$ such that 
\[\begin{cases}
f^0(s(x))=s(\P f(x))\\
f^0(t(x))=t(\P f(x))\\
\P f(x*y)=\P f(x)*\P f(y).
\end{cases}\]
The corresponding category is denoted by $\dtop$. Let \[\P_{\alpha,\beta}X = \{x\in \P X\mid s(x)=\alpha \hbox{ and } t(x)=\beta\}.\]
\ed

\begin{nota}
	The map $\P f:\P X\longrightarrow \P Y$ can be denoted by $f:\P X\to \P Y$ is there is no ambiguity. The set map $f^0:X^0\longrightarrow Y^0$ can be denoted by $f:X^0\longrightarrow Y^0$ is there is no ambiguity.
\end{nota}

\begin{exa}
	Every set can be viewed as a flow with an empty space of execution paths.
\end{exa}

One another example of flow is important for the sequel: 

\begin{exa}
	For a topological space $Z$, let $\glob(Z)$ be the flow defined by 
	\[
	\begin{cases}
	\glob(Z)^0=\{0,1\}\\
	\P \glob(Z)=\P_{0,1}\glob(Z)=Z.
	\end{cases}
	\]
This flow has no composition law.
\end{exa}

The category $\dtop$ is equipped with its \textit{q-model structure}. Its existence is proved in \cite[Theorem~7.4]{QHMmodel}. The latter paper is written in $\top_{\mathbf{\De}}$ but this result is still valid in $\top$ since the q-model structure is obtained by right-inducing a cofibrantly generated model structure using the Quillen Path Object argument. The \textit{q-model structure} of flows is the cofibrantly generated model structure such that the generating cofibrations are the maps of the form $\glob(\mathbf{S}^{n-1})\subset \glob(\mathbf{D}^{n})$ for $n\geq 0$ and the maps $C:\varnothing\subset \{0\}$ and $R:\{0,1\}\to \{0\}$, such that the weak equivalences are the maps of flows $f:X\to Y$  inducing a bijection $f^0:X^0\iso Y^0$ and a weak homotopy equivalence $\P f:\P X \to \P Y$, and such that the fibrations are the maps of flows $f:X\to Y$  inducing a q-fibration $\P f:\P X \to \P Y$. This model structure is left determined \cite[Theorem~4.3]{leftdetflow}.

This section is devoted to calculating the space of execution paths of the pushout of a flow along a map of the form $\glob(\de Z) \subset \glob(Z)$ where the map $\de Z \to Z$ is any continuous map. It is not even assumed that the map $\de Z \to Z$ is one-to-one in this section. The notations are chosen only to tell the reader how the results of this section are going to be used in the sequel.

\bp \label{final-structure-revisited-flow}
Consider a colimit cocone $(X_i)\stackrel{\bullet}\to X$ of $\dtop$. Let $f_i:X_i\to X$ be the canonical maps. Then the set of execution paths of $X$ is the set of finite compositions of the form $(f_1\gamma_1) * \dots * (f_n\gamma_n)$ such that $\gamma_i$ is an execution path of $X_i$ for all $1\leq i \leq n$.
\ep

\bpf
Every execution path $\gamma_i$ of some $X_i$ gives rise to an execution path $f_i\gamma_i$ of $X$. Every execution path of $X$ can be written as a finite composition of the form $(f_1\gamma_1) * \dots * (f_n\gamma_n)$ because of the universal property satisfied by $X$. 
\epf

Let $\de Z \to Z$ be a continuous map. Consider a pushout diagram of flows 
\[
\xymatrix@C=4em@R=4em
{
	\glob(\partial Z) \fd{} \fr{g} & A \ar@{->}[d]^-{f} \\
	\glob(Z) \fr{\widehat{g}} & \cocartesien X.
}
\]
Let $T$ be the topological space defined by the pushout diagram of $\top$
\[
\xymatrix@C=4em@R=4em
{
	\partial Z  \fd{} \fr{g} & \P_{g(0),g(1)} A \ar@{->}[d]^-{f} \\
	Z  \fr{\widehat{g}} & \cocartesien T.
}
\]
Consider the diagram of spaces $\D^f:\mathcal{P}^{g(0),g(1)}(A^0)\to \top$ defined as follows:
\[
\D^f((u_0,\epsilon_1,u_1),(u_1,\epsilon_2,u_2),\dots ,(u_{n-1},\epsilon_n,u_n)) = Z_{u_0,u_1}\p Z_{u_1,u_2} \p \dots \p Z_{u_{n-1},u_n}
\]
with 
\[
Z_{u_{i-1},u_i}=
\begin{cases}
\P_{u_{i-1},u_i}A & \hbox{if }\epsilon_i=0\\
T & \hbox{if }\epsilon_i=1
\end{cases}
\] 
In the case $\epsilon_i=1$, $(u_{i-1},u_i)=(g(0),g(1))$ by definition of $\mathcal{P}^{g(0),g(1)}(A^0)$. The inclusion maps $I_i's$ are induced by the map $f:\P_{g(0),g(1)} A \to T$. The composition maps $c_i's$ are induced by the compositions of paths of $A$. 

\bp \label{well-defined-diag}
We obtain a well-defined diagram of spaces $\D^f:\mathcal{P}^{g(0),g(1)}(A^0)\to \top$.
\ep

\bpf The relations $c_ic_j=c_{j-1}c_i$ for $i+1<j$ and $I_iI_j=I_jI_i$ for $i\neq j$ and \[c_iI_j = \begin{cases}
I_{j-1}c_i&\hbox{if } j\geq i+2\\
I_jc_i&\hbox{if } j\leq i-1
\end{cases}\] are straightforward. The relations $c_i.c_{i+1}=c_{i}.c_i$ come from the associativity of the composition of paths of $A$. 
\epf

Let $\mathcal{P}_{\alpha,\beta}^{g(0),g(1)}(A^0)\subset \mathcal{P}^{g(0),g(1)}(A^0)$ be the full subcategory generated by the objects \[((u_0,\epsilon_1,u_1),(u_1,\epsilon_2,u_2),\dots ,(u_{n-1},\epsilon_n,u_n))\] such that $u_0=\alpha$ and $u_n=\beta$ for $(\alpha,\beta)\in A^0\p A^0$. For $(\alpha,\beta)\in A^0\p A^0$, the inclusion functor $\mathcal{P}_{\alpha,\beta}^{g(0),g(1)}(A^0) \subset \mathcal{P}^{g(0),g(1)}(A^0)$ induces a well-defined diagram \[\D_{\alpha,\beta}^f:\mathcal{P}_{\alpha,\beta}^{g(0),g(1)}(A^0) \subset \mathcal{P}^{g(0),g(1)}(A^0)\longrightarrow \top.\]  We obtain a map in $\D \top$ (see Appendix~\ref{alldg})
\[
\big(\mathcal{P}_{\alpha,\beta}^{g(0),g(1)}(A^0) \subset \mathcal{P}^{g(0),g(1)}(A^0),\id\big):\D_{\alpha,\beta}^f\longrightarrow \D^f.
\]
By the universal property of the sum, we obtain a map 
\[
\bigsqcup_{(\alpha,\beta)\in A^0\p A^0} \D_{\alpha,\beta}^f \stackrel{\iso}\longrightarrow D^f
\]
which is an isomorphism by Proposition~\ref{plusdk} using the decomposition in $\dcat$ 
\[
\mathcal{P}^{g(0),g(1)}(A^0) = \bigsqcup_{(\alpha,\beta)\in A^0\p A^0} \mathcal{P}_{\alpha,\beta}^{g(0),g(1)}(A^0).
\]

\bp \label{composition_diagram}
For all triples $(\alpha,\beta,\gamma)\in A^0\p A^0\p A^0$, the concatenation of tuples induces a functor
$*:\mathcal{P}_{\alpha,\beta}^{g(0),g(1)}(A^0) \p \mathcal{P}_{\beta,\gamma}^{g(0),g(1)}(A^0)\longrightarrow \mathcal{P}_{\alpha,\gamma}^{g(0),g(1)}(A^0)$.
\ep

\bpf The small categories $\mathcal{P}_{\alpha,\beta}^{g(0),g(1)}(A^0)$, $\mathcal{P}_{\beta,\gamma}^{g(0),g(1)}(A^0)$ and $\mathcal{P}_{\alpha,\gamma}^{g(0),g(1)}(A^0)$ are posets by Proposition~\ref{poset}. Let $\underline{m},\underline{n}\in \mathcal{P}_{\alpha,\beta}^{g(0),g(1)}(A^0)$ and $\underline{m'},\underline{n'}\in \mathcal{P}_{\beta,\gamma}^{g(0),g(1)}(A^0)$ with $\underline{m}\leq \underline{n}$ and $\underline{m'}\leq \underline{n'}$. Using the relations of Group C, let 
\begin{align*}
& \underline{n}=(I_{j_1}\dots I_{j_q})(c_{i_1}\dots c_{i_p})(\underline{m}), \\
& \underline{n'}=(I_{j'_1}\dots I_{j'_q})(c_{i'_1}\dots c_{i'_p})(\underline{m'}).
\end{align*}
Then there is the equality
\[
(\underline{n},\underline{n'}) = (I_{j'_1+w}\dots I_{j'_q+w})(c_{i'_1+w}\dots c_{i'_p+w})(I_{j_1}\dots I_{j_q})(c_{i_1}\dots c_{i_p})(\underline{m},\underline{m'})
\]
where $w$ is the length of $\underline{n}$. We deduce that $(\underline{m},\underline{m'})\leq (\underline{n},\underline{n'})$ in $\mathcal{P}_{\alpha,\gamma}^{g(0),g(1)}(A^0)$.
\epf

We therefore obtain a map in $\D \top$ 
\[
\big(*:\mathcal{P}_{\alpha,\beta}^{g(0),g(1)}(A^0) \p \mathcal{P}_{\beta,\gamma}^{g(0),g(1)}(A^0)\longrightarrow \mathcal{P}_{\alpha,\gamma}^{g(0),g(1)}(A^0),\id\big):\D_{\alpha,\beta}^f \p \D_{\beta,\gamma}^f \longrightarrow \D_{\alpha,\gamma}^f
\]
for all $(\alpha,\beta,\gamma)\in A^0 \p A^0 \p A^0$. Using Proposition~\ref{cartesiandk}, we obtain a continuous map \[*:\liminj \D_{\alpha,\beta}^f \p \liminj \D_{\beta,\gamma}^f \longrightarrow \liminj \D_{\alpha,\gamma}^f.\]
Since the concatenation of tuples is associative, we obtain a well-defined flow $\overline{X}$ by setting $\overline{X}^0=A^0$, $\P_{\alpha,\beta}\overline{X} = \liminj \D_{\alpha,\beta}^f$ and with the composition law above.

\bth \label{calculation-pathspace-flow} (replacement for \cite[Proposition~15.1]{model3}) With the notations above. We obtain a commutative square of maps of flows 
\[\xymatrix@C=4em@R=4em
{
	\glob(\partial Z) \fd{} \fr{} & A \ar@{->}[d]^-{} \\
	\glob(Z) \fr{} &  \overline{X}
}\]
which is a pushout diagram. In particular, we obtain the homeomorphism \[\liminj \D^f \iso \P X.\]
\eth

\bpf
The map $\glob(Z)\to \overline{X}$ induced by the mapping $z\in Z\mapsto \widehat{g}(z)\in Z_{g(0),g(1)}$ is necessarily a map of flows since the globe does not contain composable execution paths. The map $A\to \overline{X}$ induced by the identities of $\P_{\alpha,\beta}A$ for all $(\alpha,\beta)\in A^0\p A^0$ is a map of flows because of the presence of the composition maps in $\mathcal{P}^{g(0),g(1)}(A^0)$. The square of maps of the statement of the theorem is commutative because of the presence of the inclusion maps in $\mathcal{P}^{g(0),g(1)}(A^0)$. Consider another commutative diagram of flows 
\[
\xymatrix@C=4em@R=4em
{
	\glob(\partial Z) \fd{} \fr{g} & A \ar@{->}[d]^-{h} \\
	\glob(Z) \fr{k} &  U.
}
\]
It induces the commutative diagram of topological spaces
	\[
\xymatrix@C=4em@R=4em
{
	\partial Z  \fd{} \fr{g} & \P_{g(0),g(1)} A \ar@/^20pt/@{->}[rrdd]^-{\P h}\ar@{->}[d]^-{f} \\
	Z  \fr{\widehat{g}} \ar@/_20pt/@{->}[rrrd]_-{\P k} & \cocartesien T\ar@{-->}[rrd]|-{\psi}&  \circled{A}\\
	&&& \P_{h(0),h(1)} U.
}
\]
The universal property of the pushout yields a map $\psi:T \to \P_{g(0),g(1)}U$. We obtain a composite continuous map 
\begin{multline*}
\D^f((u_0,\epsilon_1,u_1),(u_1,\epsilon_2,u_2),\dots ,(u_{n-1},\epsilon_n,u_n))
\\ \longrightarrow \P_{u_{0},u_1}U \p \dots \p \P_{u_{n-1},u_{n}}U  \longrightarrow \P_{u_{0},u_n}U,
\end{multline*}
where the left-hand map is a product of $\psi{'s}$ and $\P h{'s}$ and where the right-hand map is the composition law of $U$. Thanks to the naturality of the composition, and thanks to the commutativity of the triangle $\circled{A}$, we obtain a cocone \[\D^f \stackrel{\bullet}\longrightarrow \P U\] and therefore a canonical map $\liminj \D^f \longrightarrow \P U$. It is straightforward to verify that we have obtained a well-defined map of flows from $\widehat{\psi}:\overline{X}\to U$ such that the following diagram of flows is commutative:
\[
\xymatrix@C=4em@R=4em
{
	\glob(\partial Z)  \fd{} \fr{g} &  A \ar@/^20pt/@{->}[rrdd]^-{h}\ar@{->}[d]^-{} \\
	\glob(Z)  \fr{} \ar@/_20pt/@{->}[rrrd]_-{k} &  \overline{X} \ar@{-->}[rrd]|-{\widehat{\psi}}&  \\
	&&&  U.
}
\]
The uniqueness of $\widehat{\psi}:\overline{X}\to U$ is a consequence of Proposition~\ref{final-structure-revisited-flow}. 
\epf

\section{Left properness}
\label{proof-leftproper-section}

We recall the explicit calculation of the pushout product of several morphisms. 

\bp\label{calculpushout} \cite[Theorem~B.3]{3eme} 
Let $f_i:A_i\longrightarrow B_i$ for $0\leq i\leq p$ be $p+1$ morphisms
of a bicomplete cartesian closed category $\C$. Let $S\subset
\{0,\dots,p\}$. Let 
\[C_p(S):=\lp \prod_{i\in S} B_i\rp \p \lp \prod_{i\notin S} A_i\rp.\]
If $S$ and $T$ are two subsets of $\{0,\dots,p\}$ such that $S\subset
T$, let \[C_p(i_S^T):C_p(S)\longrightarrow C_p(T)\] be the morphism 
\[\lp\prod_{i\in S}\id_{B_i}\rp\p \lp\prod_{i\in T\backslash S} f_i\rp\p 
\lp\prod_{i\notin T} \id_{A_i}\rp. \]
Then:  
\begin{enumerate} 
	\item the mappings $S\mapsto C_p(S)$ and $i_S^T\mapsto C_p(i_S^T)$ 
	give rise to a functor from the order
	complex of the poset $\{0,\dots,p\}$ to $\C$
	\item there exists a canonical morphism 
	{
		\[\liminj_{S\subsetneqq \{0,\dots,p\}} 
		C_p(S)\longrightarrow C_p(\{0,\dots,p\}).\]}
	and it is equal to the morphism $f_0\square\dots \square f_p$. 
\end{enumerate}
\ep

\bp \label{prep1}
 Let $\de Z \to Z$ be a continuous map. Consider a pushout diagram of flows 
 \[
 \xymatrix@C=4em@R=4em
 {
 	\glob(\partial Z) \fd{} \fr{g} & A \ar@{->}[d]^-{f} \\
 	\glob(Z) \fr{\widehat{g}} & \cocartesien X.
 }
 \]
 Let $T$ be the topological space defined by the pushout diagram of $\top$
 \[
 \xymatrix@C=4em@R=4em
 {
 	\partial Z \fd{} \fr{g} & \P_{g(0),g(1)} A \ar@{->}[d]^-{f} \\
 	Z \fr{\widehat{g}} & \cocartesien T.
 }
 \]
 Let $\D^f:\mathcal{P}^{g(0),g(1)}(A^0)\to \top$ be the diagram of spaces defined above: 
 \begin{itemize}[leftmargin=*]
 	\item $
 	\D^f((u_0,\epsilon_1,u_1),(u_1,\epsilon_2,u_2),\dots ,(u_{n-1},\epsilon_n,u_n)) = Z_{u_0,u_1}\p Z_{u_1,u_2} \p \dots \p Z_{u_{n-1},u_n}
 	$
 	with 
 	\[
 	Z_{u_{i-1},u_i}=
 	\begin{cases}
 	\P_{u_{i-1},u_i}A & \hbox{if }\epsilon_i=0\\
 	T & \hbox{if }\epsilon_i=1 \hbox{ (in this case, $(u_{i-1},u_i)=(g(0),g(1))$).}
 	\end{cases}
 	\] 
 	\item The composition maps $c_i's$ are induced by the compositions of paths of $A$.
 	\item The inclusion maps $I_i's$ are induced by the map $f:\P_{g(0),g(1)} A \to T$. 
 \end{itemize}
Let $\underline{n}\in \Obj(\mathcal{P}^{g(0),g(1)}(A^0))$ with $\underline{n} = ((u_0,\epsilon_1,u_1),(u_1,\epsilon_2,u_2),\dots ,(u_{n-1},\epsilon_n,u_n))$. Then the continuous map \[L_{\underline{n}} \D^f \longrightarrow \D^f(\underline{n})\] is the pushout product of the maps $\varnothing \to \P_{u_{i-1},u_i}A$ for $i$ running over $\{i\in [1,n]| \epsilon_i = 0\}$ and of the maps $\P_{g(0),g(1)} A \to T$ for $i$ running over $\{i\in [1,n]| \epsilon_i = 1\}$. Moreover, if for all $i\in [1,n]$, we have $\epsilon_i=0$, then $L_{\underline{n}} \D^f = \varnothing$. 
\ep

\bpf
It is a consequence of Proposition~\ref{calculpushout}. 
\epf

\bth \label{onestep} Let $\de Z \to Z$ be a continuous map. Consider a commutative diagram of flows: 
\[
\xymatrix@C=4em@R=4em
{
	\glob(\partial Z) \fd{} \fr{g} & A \ar@{->}[d]^-{f} \fr{s} &  A' \fd{f'} \\
	\glob(Z) \fr{\widehat{g}} & \cocartesien X \fr{\widehat{s}} & \cocartesien X'
}
\]
Suppose that the map $\de Z \to Z$ is a q-cofibration of $\top$ and that $s$ is a weak equivalence of the q-model structure of $\dtop$. Then $\widehat{s}$ is a  weak equivalence of the q-model structure of $\dtop$. 
\eth

\bpf
Since $s$ is a weak equivalence of $\dtop$, it induces a bijection $A^0\iso A'^0$. Thus we have the bijections of sets $A^0 \iso A'^0 \iso X^0 \iso X'^0$. Consider the following commutative diagram: 
\[
\xymatrix@C=4em@R=4em
{
	\partial Z  \fd{} \fr{\P g} & \P_{g(0),g(1)} A  \fr{\P s} \ar@{->}[d]^-{f} &  \P_{(sg)(0),(sg)(1)} A' \fd{f'}\\
	Z \fr{\widehat{g}} & \cocartesien T \fr{\widehat{sg}} & \cocartesien T'.
}
\]
Since the q-model structure of $\top$ is left proper, we deduce that the continuous map $\widehat{sg}:T \to T'$ is a weak homotopy equivalence. By Theorem~\ref{calculation-pathspace-flow}, there exist two diagrams $\D^f:\mathcal{P}^{g(0),g(1)}(A^0)\to \top$ and $\D^{f'}:\mathcal{P}^{g(0),g(1)}(A^0)\to \top$ and a map of diagrams $\mathcal{S}:\D^f\to \D^{f'}$ such that the map $\liminj \mathcal{S}:\liminj \D^{f} \to \liminj \D^{f'}$ is the map $\P\widehat{s}: \P X\to  \P X'$. Since $s$ is weak equivalence of the q-model structure of $\dtop$ by hypothesis, all maps \[\P s:\P_{g(\alpha),g(\beta)}  A  \to \P_{(sg)(\alpha),(sg)(\beta)}  A'\] for $(\alpha,\beta)$ running over $A^0\p A^0$ are weak homotopy equivalences. Since the binary product of two weak homotopy equivalences is a weak homotopy equivalence as well, we deduce that the map of diagrams $\mathcal{S}:\D^f\to \D^{f'}$ is a pointwise weak homotopy equivalence. By Proposition~\ref{prep1}, for all $\underline{n}\in \Obj(\mathcal{P}^{g(0),g(1)}(A^0))$, the map $L_{\underline{n}}\D^f \to \D^f(\underline{n})$ ($L_{\underline{n}}\D^{f'} \to \D^{f'}(\underline{n})$ resp.) is a pushout product of h-cofibrations of the form $\varnothing\to \P_{\alpha,\beta}A$ and of q-cofibrations of the form $f:\P_{g(0),g(1)} A\to T$ (of q-cofibrations of the form $f':\P_{(sg)(0),(sg)(1)} A\to T'$ resp.). We deduce that the diagrams $\D^{f}$ and $\D^{f'}$ are Reedy h-cofibrant, i.e. Reedy cofibrant for the Reedy model structure on the category of diagrams $\Dcat(\mathcal{P}^{g(0),g(1)}(A^0),\top)$ with $\top$ equipped with the h-model structure. Thus $\liminj \D^{f}$ ($\liminj \D^{f'}$ resp.) is the homotopy colimit of $\D^{f}$ (of $\D^{f'}$ resp.) calculated in the h-model structure of $\top$ by Theorem~\ref{ok}. By Theorem~\ref{hocolimfacile1}, these homotopy colimit have the same weak homotopy types as the homotopy colimit of $\D^{f}$ (of $\D^{f'}$ resp.) calculated in the q-model structure of $\top$. We deduce using the $2$-out-of-$3$ axiom that the map $\P \widehat{s}: \P X\to  \P X'$ is a weak homotopy equivalence and that $\widehat{s}:X\to X'$ is a weak equivalence of the q-model structure of $\dtop$. 
\epf

\bth \label{pretower} Let $\de Z \to Z$ be a continuous map. Consider a commutative diagram of flows: 
\[
\xymatrix@C=4em@R=4em
{
	\glob(\partial Z) \fd{} \fr{g} & A \ar@{->}[d]^-{f} \\
	\glob(Z) \fr{\widehat{g}} & \cocartesien X.
}
\]
Then we have: 
\begin{enumerate}[leftmargin=*]
	\item Suppose that the map $\de Z \to Z$ is a h-cofibration of $\top$. Then the map \[\P f:\P A \to \P X\] is a h-cofibration of topological spaces.
	\item Suppose that the map $\de Z \to Z$ is a trivial h-cofibration of $\top$. Then the map \[\P f:\P A \to \P X\] is a trivial h-cofibration of topological spaces.
	\item Suppose that the map $\de Z \to Z$ is a m-cofibration (a q-cofibration resp.) of $\top$ and that $\P A$ is a m-cofibrant space (a q-cofibrant space resp.). Then \[\P f:\P A \to \P X\] is a m-cofibration (a q-cofibration resp.) of topological spaces.
	\item Suppose that the map $\de Z \to Z$ is a trivial m-cofibration (a trivial q-cofibration resp.) of $\top$ and that $\P A$ is a m-cofibrant space (a q-cofibrant space resp.). Then \[\P f:\P A \to \P X\] is a trivial m-cofibration (a trivial q-cofibration resp.) of topological spaces.
\end{enumerate}
\eth

\bpf
With the notations above. The particular case $\de Z = Z$, $f=\id_A$ and $A=X$ yields the homeomorphism \[\liminj \D^{\id_A} \iso \P A.\] We have a map of diagrams $\D^{\id_A} \to \D^{f}$ which induces for all $\underline{n}\in \Obj(\mathcal{P}^{g(0),g(1)}(A^0))$ a continuous map \[L_{\underline{n}} \D^f \sqcup_{L_{\underline{n}} \D^{\id_A}} \D^{\id_A}(\underline{n}) \longrightarrow \D^{f}(\underline{n}).\] Let \[\underline{n} = ((u_0,\epsilon_1,u_1),(u_1,\epsilon_2,u_2),\dots ,(u_{n-1},\epsilon_n,u_n)).\] There are two mutually exclusive cases: 
\begin{enumerate}
	\item[(a)] All $\epsilon_i$ for $i=1,\dots n$ are equal to zero.
	\item[(b)] There exists $i\in [1,n]$ such that $\epsilon_i = 1$.
\end{enumerate}
In the case (a), we have \[\D^{\id_A}(\underline{n}) = \D^{f}(\underline{n}) = \P_{u_0,u_1}A \p \dots \p \P_{u_{n-1},u_n}A.\] Moreover, by Proposition~\ref{prep1}, we have $L_{\underline{n}} \D^{\id_A} = L_{\underline{n}} \D^f = \varnothing$. We deduce that the map \[L_{\underline{n}} \D^f \sqcup_{L_{\underline{n}} \D^{\id_A}} \D^{\id_A}(\underline{n}) \longrightarrow \D^{f}(\underline{n})\] is isomorphic to the identity of $\D^{\id_A}(\underline{n})$. In the case (b), The map \[L_{\underline{n}} \D^{\id_A} \longrightarrow \D^{\id_A}(\underline{n})\] is by Proposition~\ref{prep1} a pushout product of several maps such that one of them is the identity map $\id: \P _{g(0),g(1)}A \to \P _{g(0),g(1)}A$ because $\epsilon_i=1$ for some $i$. Therefore the map $L_{\underline{n}} \D^{\id_A} \to \D^{\id_A}(\underline{n})$ is an isomorphism. We deduce that the map \[L_{\underline{n}} \D^f \sqcup_{L_{\underline{n}} \D^{\id_A}} \D^{\id_A}(\underline{n}) \longrightarrow \D^{f}(\underline{n})\] is isomorphic to the map $L_{\underline{n}} \D^f \to \D^f(\underline{n})$. By Proposition~\ref{prep1}, for all objects $\underline{n}\in \Obj(\mathcal{P}^{g(0),g(1)}(A^0))$, the map $L_{\underline{n}}\D^f \to \D^f(\underline{n})$ is a pushout product of maps of the form $\varnothing\to \P _{\alpha,\beta}A$ and of the form $f:\P_{g(0),g(1)} A\to T$. We conclude that the map 
\[L_{\underline{n}} \D^f \sqcup_{L_{\underline{n}} \D^{\id_A}} \D^{\id_A}(\underline{n}) \longrightarrow \D^{f}(\underline{n})\]
is for all $\underline{n}$ either an isomorphism, or a pushout product of maps of the form $\varnothing\to \P _{\alpha,\beta}A$ and of the form $f:\P_{g(0),g(1)} A\to T$, the latter appearing at least once in the pushout product and being a pushout of the map $\de Z\subset Z$. We are now ready to complete the proof.

\textbf{Case (1)}. 
The map $L_{\underline{n}}\D^f \to \D^f(\underline{n})$ is a h-cofibration of $\top$ for all $\underline{n}$. We deduce that the map of diagrams $\D^{\id_A} \to \D^{f}$ is a Reedy h-cofibration. Therefore by passing to the colimit which is a left Quillen adjoint by Theorem~\ref{ok}, we deduce that the map $\P A \to \P X$ is a h-cofibration of $\top$. 

\textbf{Case (2)}. 
If moreover the map $\de Z \subset Z$ is a homotopy equivalence, then the map $L_{\underline{n}}\D^f \to \D^f(\underline{n})$ is always a trivial h-cofibration of $\top$ for all $\underline{n}$. We deduce that the map of diagrams $\D^{\id_A} \to \D^{f}$ is a trivial Reedy h-cofibration. Therefore by passing to the colimit which is a left Quillen adjoint by Theorem~\ref{ok}, we deduce that the map $\P A \to \P X$ is a trivial h-cofibration of $\top$.

\textbf{Case (3)}. Suppose now that $\P A$ is a m-cofibrant space. Then by \cite[Corollary~3.7]{mixed-cole}, the space $\P A$ is homotopy equivalent to a q-cofibrant space $U$. Each connected component of $U$ is q-cofibrant. Therefore for all $(\alpha,\beta)\in A^0\p A^0$, the topological space $\P_{\alpha,\beta} A$ is m-cofibrant. We deduce that the map $L_{\underline{n}}\D^f \to \D^f(\underline{n})$ is always a m-cofibration of $\top$ for all $\underline{n}$ because a pushout product of two m-cofibrations of spaces is a m-cofibration by \cite[Proposition~6.6]{mixed-cole}. We deduce that the map of diagrams $\D^{\id_A} \to \D^{f}$ is a Reedy m-cofibration. Therefore by passing to the colimit which is a left Quillen adjoint by Theorem~\ref{ok}, we deduce that the map $\P A \to \P X$ is a m-cofibration of $\top$. The proof is similar for the other case.

\textbf{Case (4)}. Assume moreover that the map $\de Z\subset Z$ is a weak homotopy equivalence. then the map $L_{\underline{n}}\D^f \to \D^f(\underline{n})$ is always a trivial m-cofibration of $\top$ for all $\underline{n}$. We deduce that the map of diagrams $\D^{\id_A} \to \D^{f}$ is a trivial Reedy m-cofibration. Therefore by passing to the colimit which is a left Quillen adjoint by Theorem~\ref{ok}, we deduce that the map $\P A \to \P X$ is a trivial m-cofibration of $\top$. The proof is similar for the other case.
\epf

\bth  \label{path-almost-accessible}
Let $X:\lambda \to \dtop$ be a transfinite tower of flows such that for all $\mu<\lambda$, the map $\P X_\mu \to \P X_{\mu+1}$ is a relative-$T_1$ inclusion. Then the canonical map \[\liminj (\P.X) \longrightarrow \P \liminj X\] is a homeomorphism. Moreover the topology of $\P \liminj X$ is the final topology.
\eth

\bpf Let us use the notation $\dtop(\top_\B)$ to specify that the underlying category is $\top_\B$. 

1) Let $Y = \liminj\nolimits^{\dtop(\top_\B)} X$, the colimit being taken in $\dtop(\top_\B)$. The forgetful functor $\omega:\top_\B\to\set$ induces a forgetful functor $\widehat{\omega}:\dtop(\top_\B)\to\dtop(\set)$~\footnote{The objects of $\dtop(\set)$ are the small semicategories, or small non-unital categories.}. In particular, we have $\P\widehat{\omega}(Z) = \omega(\P Z)$ for all $Z$, $\widehat{\omega}(Z)^0=Z^0$ and the composition law of $\widehat{\omega}(Z)$ is the composition law of $Z$. The functor \[\widehat{\omega}:\dtop(\top_\B)\to\dtop(\set)\] has a right adjoint by equipping the set of paths with the indiscrete topology (which is a $\B$-generated space): since every map to an indiscrete space is continuous, the composition law is automatically continuous. Using Proposition~\ref{final-structure-revisited-flow}, we observe that the functor \[\P:\dtop(\set) \to \set\] is finitely accessible. 

2) We can now derive the following sequence of bijections of sets: 
\begin{align*}
\omega(\P Y) &= \P (\widehat{\omega}Y) & \hbox{ by definition of $\widehat{\omega}$}\\
&\iso \P \liminj(\widehat{\omega}.X) & \hbox{ since $\widehat{\omega}$ is colimit-preserving}\\
&\iso  \liminj(\P.\widehat{\omega}.X) & \hbox{ since $\P:\dtop(\set) \to \set$ is finitely accessible}\\
&\iso \liminj(\omega.\P.X) & \hbox{ by definition of $\widehat{\omega}$}\\
&\iso \omega\big(\liminj\nolimits^{\top_\B}(\P.X)\big) & \hbox{ since ${\omega}$ is colimit-preserving.}
\end{align*}
In other terms, the two topological spaces $\P Y$ and $\liminj\nolimits^{\top_\B}(\P.X)$ have the same underlying set. We deduce that the canonical map $\liminj^{\top_\B} (\P.X) \to \P Y$ is a continuous bijection. The left-hand space is equipped with the final topology (the symbol $\liminj^{\top_\B}$ means that the colimit is calculated in $\top_\B$). The composite law of $\widehat{\omega}Y$ is defined as follows. Let $x,y\in \omega \P Y$ be two composable paths. Since the colimit is directed, it is possible to find $x_i,y_i\in \omega \P X_i$ for some $i<\lambda$ composable such that the canonical map $\P X_i \to \P Y$ takes $x_i$ to $x$ and $y_i$ to $y$. Then we set $x*y$ to be the image of $x_i*y_i$ by the canonical map $\P X_i \to \P Y$. It is a standard argument to prove that this yields a well-defined associative composition law on $\omega \P Y$. 

3) The next step is to prove that the above set map $\omega\P Y \p_{Y^0} \omega\P Y \to \omega\P Y$ is continuous if $\omega\P Y$ is equipped with the final topology. Let $K\in \B$.  Since we are working with $\B$-generated spaces, it suffices to prove that for any continuous map $K \to \P Y \p_{Y^0} \P Y$, the composite map $K \to \P Y \p_{Y^0} \P Y \to \P Y$ is continuous. First note that for every limit ordinal $\mu\leq \lambda$, $\P X_\mu$ is always equipped with the final topology because all maps $\P X_\mu \to \P X_{\mu+1}$ for $\mu<\lambda$ are one-to-one by hypothesis. Since all maps $\P X_\mu \to \P X_{\mu+1}$ for $\mu<\lambda$ are relative-$T_1$ by hypothesis, there exists an ordinal $\nu<\lambda$ and a commutative diagram of $\B$-generated spaces
\[
\xymatrix@C=4em@R=4em
{
K \fr{} \ar@{=}[d] & \P X_\nu \p_{X_\nu^0} \P X_\nu  \fr{*} \fd{} & \P X_\nu \fd{}\\
K \fr{} & \P Y \p_{Y^0} \P Y  \fr{*} & \P Y
}
\]
The top arrow $\P X_\nu \p_{X_\nu^0} \P X_\nu \to \P X_\nu$ is continuous because it is the composition law of a flow. We deduce that the bottom arrow $\P Y \p_{Y^0} \P Y  \to \P Y$ is continuous as well by equipping $\P Y$ with the final topology.

4) The flow $Y = \liminj\nolimits^{\dtop(\top_\B)} X$ satisfies the universal property of the colimit in $\dtop$. The universal property of the colimit yields a map of flows $Y \to \liminj X$ such that the composite $Y \to \liminj X \to Y$ is the identity of $Y$. Thus, the canonical map $\liminj^{\top_\B} (\P.X) \to \P Y$ is a homeomorphism.
\epf

\bth \label{proof-left-proper-flow}
The q-model category of flows is left proper.
\eth

\begin{figure}
	\[
\xymatrix@C=5em
{
	A_0=A \fr{s_0=s} \fd{} & A'_0=A' \fd{} \\
	\vdots \fd{} & \vdots \fd{} \\
	A_{\lambda}\fr{s_{\lambda}}\fd{}& A'_{\lambda} \fd{} \\
	A_{\lambda+1}\fr{s_{\lambda+1}}\fd{}& \cocartesien A'_{\lambda+1} \fd{}\\
	\vdots\fd{} & \vdots\fd{} \\
	\liminj A_{\lambda} = X \fr{\liminj s_{\lambda}} & \liminj A'_{\lambda} = X'
}
\]
\caption{}
\label{tower1}
\end{figure}

\begin{figure}
	\[
\xymatrix@C=5em
{
	\P A_0 \fr{\P  s_0} \fd{} & \P A'_0 \fd{} \\
	\vdots \fd{} & \vdots \fd{} \\
	\P A_{\lambda}\fr{\P s_{\lambda}}\fd{}& \P A'_{\lambda} \fd{} \\
	\P A_{\lambda+1}\fr{\P s_{\lambda+1}}\fd{}& \P A'_{\lambda+1} \fd{}\\
	\vdots\fd{} & \vdots\fd{} \\
	\liminj \P A_{\lambda} = \P X \fr{\liminj \P  s_{\lambda}} & \liminj \P  A'_{\lambda} = \P X'.
}
\]
\caption{}
\label{tower2}
\end{figure}

\bpf Consider the commutative diagram
\[
\xymatrix@C=4em@R=4em
{
	U \fd{i} \fr{g} & A \ar@{->}[d]^-{f} \fr{s} &  A' \fd{f'} \\
	V \fr{\widehat{g}} & \cocartesien X \fr{\widehat{s}} & \cocartesien X'
}
\]
where $i:U\to V$ is a transfinite composition of q-cofibrations of the form 
\begin{enumerate}
	\item $C:\varnothing \subset \{0\}$,
	\item $R:\{0,1\} \to \{0\}$, 
	\item $\glob(\de Z)\to \glob(Z)$ where the inclusion $\de Z\subset Z$ is a q-cofibration of spaces.
\end{enumerate}
We obtain a map of transfinite towers of flows as depicted in Figure~\ref{tower1}. Each map $\P A_\lambda \to \P A_{\lambda+1}$ is a h-cofibration of $\top$, and therefore is relative-$T_1$ for the following reasons: 
\begin{enumerate}
	\item if $A_\lambda \to A_{\lambda+1}$ is a pushout of $C:\varnothing \subset \{0\}$, then $\P A_\lambda = \P A_{\lambda+1}$.
	\item if $A_\lambda \to A_{\lambda+1}$ is a pushout of $R:\{0,1\} \to \{0\}$, then $\P A_{\lambda+1} \iso \P A_\lambda \sqcup W$ for some topological space $W$ (the space of execution paths freely generated by identifying two states; $W$ is empty when the map $\{0,1\}\to A_\lambda$ is constant). 
	\item if $A_\lambda \to A_{\lambda+1}$ is a pushout of $\glob(\de Z)\to \glob(Z)$ where the inclusion $\de Z\subset Z$ is a q-cofibration, then $\P A_\lambda \to \P A_{\lambda+1}$ is a h-cofibration of $\top$ by Theorem~\ref{pretower}~(1). 
\end{enumerate}
By Theorem~\ref{path-almost-accessible}, we obtain a map of transfinite towers of topological spaces as depicted in Figure~\ref{tower2}. Let us prove by transfinite induction on $\lambda\geq 0$ that $\P s_\lambda$ is a weak homotopy equivalence. The induction hypothesis holds for $\lambda=0$ by hypothesis. If it is proved for $\lambda\geq 0$, then it holds for $\lambda+1$ 
\begin{itemize}
	\item by Theorem~\ref{onestep} for a pushout of a q-cofibration of the form $\glob(\de Z)\to \glob(Z)$
	\item because $\P s_\lambda=\P s_{\lambda+1}$ for a pushout of $C:\varnothing \subset \{0\}$
	\item because the binary product of two weak homotopy equivalences is a weak homotopy equivalence for a pushout of $R:\{0,1\} \to \{0\}$ (the argument is explained in the proof of \cite[Theorem~7.4]{2eme}.
\end{itemize} 
If $\lambda$ is a limit ordinal, and if the induction hypothesis is proved for all $\mu<\lambda$, since all vertical maps are h-cofibrations of topological spaces by Theorem~\ref{pretower}~(1), the colimits are actually homotopy colimits for the h-model structure of $\top$, and by Theorem~\ref{hocolimfacile1} homotopy colimits for the q-model structure of $\top$. We deduce using the $2$-out-of-$3$ axiom the induction hypothesis for $\lambda$ and the induction is complete. To complete the proof, it suffices to remember that every q-cofibration of flows is a retract of a map like $i$.
\epf

The section concludes with some additional information about the path space functor of flows.

\bth \label{complement-flow} Let $f:X\to Y$ be a (trivial resp.) q-cofibration between q-cofibrant flows. Then the continuous map $\P f:\P X \to \P Y$ is a (trivial resp.) q-cofibration of spaces. In particular, the path space functor $\P :\dtop\to \top$ preserves q-cofibrancy. 
\eth

\bpf
Let $X$ be a q-cofibrant flow. Then the map $\varnothing\to X$ is a retract of a transfinite composition of pushouts along the generating q-cofibrations. Using Theorem~\ref{pretower}~(3) and a transfinite induction, we deduce that $\P X$ is q-cofibrant.  
\epf

Theorem~\ref{complement-flow} is quite surprising because the end of the section proves that the path space functor is a right Quillen adjoint if $\B=\mathbf{\De}$.

\begin{lem}
	Let $Z$ be a general topological space. Let $\CC(Z)$ be the set of connected components of $Z$ equipped with the final topology coming from the quotient map \[p_Z:Z\longrightarrow \CC(Z).\] Then $\CC(Z)$ is totally disconnected. If $Z$ is $\De$-generated, then $\CC(Z)$ is discrete.
\end{lem}

\bpf
Let $T\subset \CC(Z)$ be a subset containing at least two points. Then $p_Z^{-1}(T)$ is not connected. Therefore there exists a nonconstant continous map $p_Z^{-1}(T) \to \{0,1\}$ with $\{0,1\}$ equipped with the discrete topology. This map factors (uniquely) as a composite $p_Z^{-1}(T) \to T \to \{0,1\}$ since $T$ is equipped with the relative topology with respect to the final topology. We deduce that the right-hand map $T \to \{0,1\}$ is nonconstant and that $T$ is not connected. We have proved that $\CC(Z)$ is totally disconnected. Assume that $Z$ is $\De$-generated. Then, by \cite[Proposition~2.8]{mdtop}, $Z$ is homeomorphic to the disjoint sum of its nonempty connected components. Therefore the connected components are open. So every point of $\CC(Z)$ is open in $\CC(Z)$ because the latter is equipped with the final topology. We deduce that $\CC(Z)$ is discrete.
\epf

\bth \label{rightadjointpath-flow} Assume that $\top$ is either $\Htop_{\mathbf{\De}}$ or $\top_{\mathbf{\De}}$, i.e. $\B=\mathbf{\De}$. The path space functor $\P :\dtop\to \top$ is a right Quillen adjoint with $\top$ equipped with its q-model structure. In particular, the functor $\P :\dtop\to \top$ is accessible.  
\eth

\bpf That the path space functor $\P :\dtop\to \top$ takes fibrations (trivial fibrations resp.) of $\dtop$ to fibrations (trivial fibrations resp.) of $\top$ comes from the construction of the q-model structure of flows \cite[Theorem~7.4]{QHMmodel}. The left adjoint $\mathbf{G}:\top \to {\dtop}$ is defined on objects as follows: 
\begin{itemize}
	\item The space of paths $\P \mathbf{G}(Z)$ is $Z$.
	\item The discrete space of states $\mathbf{G}(Z)^0$ is $\CC(Z)\p \{0,1\}$. 
	\item The source map $s:Z\to \mathbf{G}(Z)^0$ is the composite map $Z\stackrel{p_Z}\to \CC(Z)\p\{0\} \subset \mathbf{G}(Z)^0$.
	\item The target map $t:Z\to \mathbf{G}(Z)^0$ is the composite map $Z\stackrel{p_Z}\to \CC(Z)\p\{1\} \subset \mathbf{G}(Z)^0$.
	\item There is no composition law.
\end{itemize}
The definition of $\mathbf{G}:\top \to {\dtop}$ on maps is clear. Choosing a map of flows from $\mathbf{G}(Z)$ to a flow $U$ is equivalent to choosing a continuous map from $Z$ to $\P U$ because the image of any element of $\mathbf{G}(Z)^0$ is forced. \epf

\bth Assume that $\top$ is either $\top_\mathbf{K}$ or $\Htop_\mathbf{K}$, i.e. $\B=\mathbf{K}$. The path space functor $\P :\dtop\to \top$ is not a right adjoint. 
\eth

\bpf Assume that the left adjoint $\mathbf{G}:\top\to \dtop$ exists. Let $Z$ be an arbitrary object of $\top$. By naturality of the maps, we have a commutative diagram of spaces 
\[
\xymatrix@C=5em
{
	\{0\} \fr{\eta_{\{0\}}} \fd{} & \P(\mathbf{G}(\{0\})) = \P(\mathbf{G}(\P \glob(\{0\}))) \fr{\P(\epsilon_{\glob(\{0\})})}& \P\glob(\{0\}) = \{0\} \fd{} \\
	Z \fr{\eta_Z} & \P(\mathbf{G}(Z)) = \P(\mathbf{G}(\P \glob(Z))) \fr{\P(\epsilon_{\glob(Z)})}& \P\glob(Z) = Z
}
\]
for any map $\{0\}\to Z$ where $\eta:\id\Rightarrow \P.\mathbf{G}$ is the unit of the adjunction and where $\epsilon:\mathbf{G}.\P\Rightarrow \id$ is the counit of the adjunction. Thus the bottom composite map is the identity of $Z$. We obtain the homeomorphism \[Z \iso \{z\in \P(\mathbf{G}(Z)) \mid \eta_Z(\P(\epsilon_{\glob(Z)})(z))=z\}\] because both spaces are the equalizer of the pair of maps $(\eta_Z.\P(\epsilon_{\glob(Z)}),\id_{\P(\mathbf{G}(Z))})$. Consequently, $Z$ can be identified with a subset of $\P(\mathbf{G}(Z))$ equipped with the relative topology. From the existence of the map of flows $\epsilon_{\glob(Z)}:\mathbf{G}(Z) \to \glob(Z)$, we deduce that $\mathbf{G}(Z)$ has no composition law. Indeed, if $\gamma_1$ and $\gamma_2$ are two execution paths of $\mathbf{G}(Z)$, then $t(\epsilon_{\glob(Z)}(\gamma_1))=1$ and $s(\epsilon_{\glob(Z)}(\gamma_2))=0$, and therefore $t(\gamma_1)\neq s(\gamma_2)$. We deduce the existence of a flow $U_Z \subset \mathbf{G}(Z)$ without composition law defined as follows: $\P U_Z = Z$, $U_Z^0 = s(Z) \sqcup t(Z)$ with the source map $s\rest_Z:Z\subset \P(\mathbf{G}(Z)) \stackrel{s} \to s(Z)$ and the target map $t\rest_Z:Z\subset \P(\mathbf{G}(Z)) \stackrel{t} \to t(Z)$. Since $U_Z$ has no composition law, choosing a map of flows from $U_Z$ to some flow $U$ is equivalent to choosing a continuous map from $Z=\P U_Z$ to $\P U$. Therefore, $U_Z$ satisfies the same universal property as $\mathbf{G}(Z)$. We deduce the isomorphism $U_Z\iso \mathbf{G}(Z)$ and that the inclusion $Z \subset \P \mathbf{G}(Z)$ is an equality. We are now ready to reach the contradiction. 

Let $Z=\{1/n\mid n\geq 1\}\cup \{0\}$ equipped with the relative topology induced by the one of $[0,1]$. Since $Z$ is a closed subset of $[0,1]$, it belongs to $\top_\mathbf{K}$. Since it is Hausdorff, it also belongs to $\Htop_\mathbf{K}$~\footnote{Note that $Z$ does not belong to $\top_{\mathbf{\De}}$ (see the proof of Proposition~\ref{HLocPre})}. Consider the continuous map $s:Z = \P \mathbf{G}(Z)\to \mathbf{G}(Z)^0$. Then $s^{-1}(s(0))$ is an open subset of $Z$ because $\mathbf{G}(Z)^0$ is discrete. Therefore it contains $Z^{\geq p}=\{1/n\mid n\geq p\}\cup \{0\}$ for some $p\geq 1$. This implies that $s(Z) = s(\{1,1/2,\dots,1/p\})$ is finite. On the other hand, there are the homeomorphisms $Z =  \{1\} \sqcup Z^{\geq 2} \iso \{1\} \sqcup Z$. Since $\mathbf{G}$ is a left adjoint, we have $\mathbf{G}(Z) \iso \mathbf{G}(\{1\}) \sqcup \mathbf{G}(Z)$. A map $f$ from $\glob(\{1\})$ to some flow $U$ is characterized by the choice of $f(1)\in \P U$ (the images of $0=s(1)$ and $1=t(1)$ are forced), which means that $\mathbf{G}(\{1\})\iso \glob(\{1\})$. We obtain the isomorphism $\mathbf{G}(Z) \iso \glob(\{1\})) \sqcup \mathbf{G}(Z)$. We deduce the bijection of sets $s(Z) \iso \{0\} \sqcup s(Z)$, which implies that $s(Z)$ is infinite. Contradiction. 
\epf

The path space functor $\P:\dtop\to \top$ with $\B=\mathbf{\De}$ from the q-model structure of flows to the q-model structure of topological spaces is a right Quillen adjoint which preserves q-cofibrations and trivial q-cofibrations between q-cofibrant objects. Another example of such a phenomenon has been given in MathOverflow by Simon Henry \cite{rare}: see \cite[Theorem~3.2(5)]{examplerare}.

Morally speaking, in the case $\B=\mathbf{K}$, the left adjoint $\mathbf{G}:\top\to \dtop$ exists but the functor $\mathbf{G}$ takes a space $Z$ which is not homeomorphic to the disjoint sum of its connected components to a flow outside the category $\dtop$. Let us formalize this fact in the last theorem of the section. 

\bd A {\rm generalized flow} $X$ consists of a topological space $\P X$ of execution paths, a totally disconnected space $X^0$ of states, two continuous maps $s$ and $t$ from $\P X$ to $X^0$ called the source and target map respectively, and a continuous and associative map \[*:\{(x,y)\in \P X\p \P X; t(x)=s(y)\}\longrightarrow \P X\] such that $s(x*y)=s(x)$ and $t(x*y)=t(y)$.  A morphism of generalized flows $f:X\longrightarrow Y$ consists of a continuous map $f^0:X^0\longrightarrow Y^0$ together with a continuous map $\P f:\P X\longrightarrow \P Y$ such that $f(s(x))=s(f(x))$, $f(t(x))=t(f(x))$ and $f(x*y)=f(x)*f(y)$. The corresponding category is denoted by $\overline{\dtop}$. Let $\P_{\alpha,\beta}X = \{x\in \P X\mid s(x)=\alpha \hbox{ and } t(x)=\beta\}$.
\ed

\bth The path space functor $\P:\overline{\dtop} \to \top$ is a right adjoint with $\B=\mathbf{\De}$ and with $\B=\mathbf{K}$. \eth

\bpf
The left adjoint $\mathbf{G}:\top \to \overline{\dtop}$ is a kind of generalized globe. It is defined on objects as follows: 
\begin{itemize}
	\item The space of paths $\P \mathbf{G}(Z)$ is $Z$.
	\item The totally disconnected space of states $\mathbf{G}(Z)^0$ is $\CC(Z)\p \{0,1\}$.
	\item The source map $s:Z\to \mathbf{G}(Z)^0$ is the composite map $Z\stackrel{p_Z}\to \CC(Z)\p\{0\} \subset \mathbf{G}(Z)^0$.
	\item The target map $t:Z\to \mathbf{G}(Z)^0$ is the composite map $Z\stackrel{p_Z}\to \CC(Z)\p\{1\} \subset \mathbf{G}(Z)^0$.
	\item There is no composition law.
\end{itemize}
The rest of the proof is similar to the proof of Theorem~\ref{rightadjointpath-flow}. 
\epf

\section{Erratum}
\label{flaw}

This section concludes the paper by carefully describing the flaw in the proof of \cite[Proposition~15.1]{model3} and by fixing some proofs published in past papers. The section starts by a short reminder about multipointed $d$-spaces.

\subsection*{Multipointed $d$-space}

\bd \label{multid} A {\rm multipointed space} is a pair $(|X|,X^0)$ where
\begin{itemize}
	\item $|X|$ is a topological space called the {\rm underlying space} of $X$.
	\item $X^0$ is a subset of $|X|$ called the {\rm set of states} of $X$.
\end{itemize}
A morphism of multipointed spaces $f:X=(|X|,X^0) \rightarrow Y=(|Y|,Y^0)$ is a commutative square
\[
\xymatrix{
	X^0 \fr{f^0}\fd{} & Y^0 \fd{} \\ 
	|X| \fr{|f|} & |Y|.}
\] 
The corresponding category is denoted by $\mtop$.  \ed

We have the well-known proposition: 

\bp (The Moore composition) Let $U$ be a topological space. Let $\ell_1,\dots,\ell_n\in ]0,+\infty[$ be real numbers. Let $\gamma_i:[0,\ell_i]\to U$ be $n$ continuous maps with $1\leq i \leq n$. Suppose that $\gamma_i(\ell_{i})=\gamma_{i+1}(0)$ for $1\leq i < n$ (there is nothing to verify for $n=1$). Then there exists a unique continuous map $\gamma_1*\dots *\gamma_n:[0,\sum_i\ell_i]\to U$ such that 
\[
(\gamma_1*\dots *\gamma_n)(t) = \gamma_i\big(t-\sum_{j<i}\ell_i\big) \hbox{ for }\sum_{j<i} \ell_i\leq t \leq \sum_{j\leq i} \ell_i.
\]
In particular, there is the strict equality $(\gamma_1* \gamma_2)* \gamma_3=\gamma_1*(\gamma_2* \gamma_3)$.
\ep

\begin{nota}
	Let $\mu_{\ell}:[0,\ell]\to [0,1]$ be the homeomorphism defined by $\mu_\ell(t) = t/\ell$. 
\end{nota}

\bd \label{composition_map} The map $\gamma_1*\gamma_2$ is called the {\rm (Moore) composition} of $\gamma_1$ and $\gamma_2$.  The composite 
\[\xymatrix@C=4em{\gamma_1 *_N \gamma_2: [0,1] \fr{(\mu_2)^{-1}}& [0,2]
	\fr{\gamma_1*\gamma_2}& U}\] 
is called the {\rm (normalized) composition}. The normalized composition being not associative, a notation like $\gamma_1 *_N \dots *_N \gamma_n$ will mean, by convention, that $*_N$ is applied from the left to the right. 
\ed

\bd \cite{mdtop} A {\rm multipointed $d$-space $X$} is a triple
$(|X|,X^0,\P^{\mathcal{G}}X)$ where
\begin{itemize}[leftmargin=*]
	\item The pair $(|X|,X^0)$ is a multipointed space.
	\item The set $\P^{\mathcal{G}}X$ is a set of continous maps from $[0,1]$ to $|X|$ called the {\rm execution paths}, satisfying the following axioms:
	\begin{itemize}
		\item For any execution path $\gamma$, one has $\gamma(0),\gamma(1)\in X^0$.
		\item Let $\gamma$ be an execution path of $X$. Then any composite $\gamma.\phi$ with $\phi\in \mathcal{G}$ is an execution path of $X$ where $\mathcal{G}$ is the group of nondecreasing homeomorphisms from $[0,1]$ to itself.
		\item Let $\gamma_1$ and $\gamma_2$ be two composable execution paths of $X$; then the normalized composition $\gamma_1 *_N \gamma_2$ is an execution path of $X$.
	\end{itemize}
\end{itemize}
A map $f:X\to Y$ of multipointed $d$-spaces is a map of multipointed spaces $|f|$ from $(|X|,X^0)$ to $(|Y|,Y^0)$ such that for any execution path $\gamma$ of $X$, the composite map $\P^{\mathcal{G}}(\gamma) := |f|\gamma$ is an execution path of $Y$.  The category of multipointed $d$-spaces is denoted by $\ptop{\mathcal{G}}$. The subset of execution paths from $\alpha$ to $\beta$ is the set of $\gamma\in\P^{\mathcal{G}} X$  such that $\gamma(0)=\alpha$ and $\gamma(1)=\beta$; it is denoted by $\P^{\mathcal{G}}_{\alpha,\beta} X$. It is equipped with the $\B$-kelleyfication of the initial topology making the inclusion $\P^{\mathcal{G}}_{\alpha,\beta} X\subset \ttop([0,1],|X|)$ continuous. Let \[\P^{\mathcal{G}}X := \bigsqcup_{(\alpha,\beta) \in X^0\p X^0} \P^{\mathcal{G}}_{\alpha,\beta} X.\] \ed

By definition, the topological space $\P^{\mathcal{G}}X$ is homeomorphic to the disjoint union of the topological spaces $\P^{\mathcal{G}}_{\alpha,\beta} X$ for $(\alpha,\beta)$ running over $X^0\p X^0$.

The following examples play an important role in the sequel. 
\begin{enumerate}[leftmargin=*]
	\item Any set $E$ will be identified with the multipointed $d$-space $(E,E,\varnothing)$.
	\item The \textit{topological globe of $Z$}, which is denoted by $\glob^{\mathcal{G}}(Z)$, is the multipointed $d$-space defined as follows
	\begin{itemize}
		\item the underlying topological space is the quotient space \[\frac{\{0,1\}\sqcup (Z\p[0,1])}{(z,0)=(z',0)=0,(z,1)=(z',1)=1}\]
		\item the set of states is $\{0,1\}$
		\item the set of execution paths is the set of continuous maps \[\{\delta_z\phi\mid \phi\in \mathcal{G},z\in  Z\}\]
		with $\delta_z(t) = (z,t)$. 
	\end{itemize}
	In particular, $\glob^{\mathcal{G}}(\varnothing)$ is the multipointed $d$-space $\{0,1\} = (\{0,1\},\{0,1\},\varnothing)$.  
\end{enumerate}

The \textit{q-model structure} of $\ptop{\mathcal{G}}$, constructed in \cite[Theorem~6.14]{QHMmodel} (the latter paper is written in $\top_{\mathbf{\De}}$ but this result is still valid in $\top$ since the q-model structure is obtained by right-inducing a cofibrantly generated model structure using the Quillen Path Object argument), is the cofibrantly generated model structure such that the generating cofibrations are the maps of the form $\globP(\mathbf{S}^{n-1})\subset \globP(\mathbf{D}^{n})$ for $n\geq 0$ and the maps $C:\varnothing\subset \{0\}$ and $R:\{0,1\}\to \{0\}$, such that the weak equivalences are the maps of multipointed $d$-spaces $f:X\to Y$  inducing a bijection $f^0:X^0\iso Y^0$ and a weak homotopy equivalence $\P^{\mathcal{G}}f:\P^{\mathcal{G}}X \to \P^{\mathcal{G}}Y$, and such that the fibrations are the maps of multipointed $d$-spaces $f:X\to Y$  inducing a q-fibration $\P^{\mathcal{G}}f:\P^{\mathcal{G}}X \to \P^{\mathcal{G}}Y$. This model structure is also left determined by a proof similar to the proof of \cite[Theorem~4.3]{leftdetflow}.

\subsection*{Description of the flaw} 

The diagram used in \cite[Proposition~15.1]{model3} to calculate $\P X$ in the pushout diagram of flows
\[
\xymatrix@C=4em@R=4em
{
	\glob(\partial Z) \fd{} \fr{g} & A \ar@{->}[d]^-{f} \\
	\glob(Z) \fr{\widehat{g}} & \cocartesien X
}
\]
is a subdiagram $\D'^f$ of $\D^f$ which is not cofinal in $\D^f$. With the terminology of the proof of Theorem~\ref{calculation-pathspace-flow}, $\D'^f$ is the subdiagram of $\D^f$ obtained by keeping in $\mathcal{P}^{g(0),g(1)}(A^0)$ the non-simplifiable objects, and the objects $\underline{m}$ of $\mathcal{P}^{g(0),g(1)}(A^0)$ such that there is one inclusion map $\underline{m}\to \underline{n}$ (and not a composite of inclusion maps) towards a non-simplifiable object $\underline{n}$. The inclusion functor, let us denote it by $\iota$, is not cofinal. The comma category $\underline{n}\ddownarrow \iota$ is always nonempty because all non-simplifiable objects belong to this subcategory of $\mathcal{P}^{g(0),g(1)}(A^0)$ but it can be verified that it is not necessarily connected. Moreover neither $\D'^f$ nor  $\D^f$ can have a cofinal subdiagram like the one depicted in \cite[p 590]{model3}. It is possible to find counterexamples indeed. Therefore the map $\P A\to \P X$ is not equal to a transfinite composition of the kind depicted in \cite[p 590]{model3}. At least, it is not possible to prove such a fact with this method. This  incorrect argument propagated in the papers \cite{model2} \cite{3eme}, \cite{4eme}, \cite{2eme}  and \cite{mdtop}. The rest of Section~\ref{flaw} is devoted to correcting this problem.

\subsection*{Erratum for \cite{model3}}

\cite[Proposition~15.1]{model3} is only used for the proof of \cite[Theorem~15.2]{model3}. The latter remains true anyway: it is even possible to conclude that the map $\P A \to \P X$ is a homotopy equivalence, and not only a weak homotopy equivalence, thanks to Theorem~\ref{pretower}(2). Note that the q-model structure of flows constructed in \cite{model3} can be more easily recovered by using Isaev's work \cite[Theorem~3.11]{leftdetflow} or by using the theory of bifibrations \cite[Theorem~7.4]{QHMmodel}.

\subsection*{Erratum for \cite{3eme}, \cite{4eme}, \cite{2eme} and \cite{mdtop}}

I do not know whether Proposition~A.1 of \cite{3eme}~\footnote{cited in \cite[Proposition~A.1]{4eme}, \cite[Proposition~7.1]{2eme} and \cite[Proposition~A.2]{mdtop}.} is true. This proposition is used to prove the left properness of the q-model category of flows and to prove that the path space functor $\P :\dtop\to \top$ preserves q-cofibrancy. All these facts are correctly proved in Section~\ref{proof-leftproper-section}. 

\subsection*{Erratum for \cite{model2}}

Not only do I not know whether the continuous map $\pi_n(\P A)\to \pi_n(\P X)$ of \cite[Theorem~V.3.4]{model2} is onto, but also this assertion is useless. The correct statement is: 

\bth \label{fix} (replacement for \cite[Theorem~V.3.4]{model2})
Let $i:\de Z \to Z$ be a q-cofibration of $\top$ inducing an isomorphism $\pi_k(\de Z) \iso \pi_k(Z)$ for all $k<n$ and for all base points. Consider a pushout diagram of flows 
\[
\xymatrix@C=4em@R=4em{
	\glob(\de Z) \fr{g}\fd{} & A \fd{f} \\
	\glob(Z) \fr{} & \cocartesien X }
\]
with $A$ a q-cofibrant flow. Then the continuous map $\P f:{\P}A\longrightarrow {\P} X$ induces an isomorphism $\pi_k(\P A) \iso \pi_k(\P X)$ for all $k<n$ and for all base points.
\eth 

\bpf We consider the homotopical localization $\mathbf{L}_{\mathbf{S}^{n} \subset \Di^{n+1}}\top$ of the q-model structure of $\top$ by the q-cofibration $\mathbf{S}^{n} \subset \Di^{n+1}$. By \cite[Proposition~1.5.2]{ref_model2}, the map $i:\de Z\to Z$ is a trivial cofibration of this homotopical localization. By Theorem~\ref{complement-flow}, the topological space $\P A$ is q-cofibrant. With the notations of the proof of Theorem~\ref{pretower}(4), we obtain that the map $L_{\underline{n}}\D^f \to \D^f(\underline{n})$ is always a trivial cofibration of $\mathbf{L}_{\mathbf{S}^{n} \subset \Di^{n+1}}\top$ for all $\underline{n}$. It means that the map of diagrams $\D^{\id_A}\to \D^f$ is a trivial Reedy cofibration of \[\Dcat(\mathcal{P}^{g(0),g(1)}(A^0),\mathbf{L}_{\mathbf{S}^{n} \subset \Di^{n+1}}\top).\] Therefore by passing to the colimit which is a left Quillen adjoint by Theorem~\ref{ok}, we deduce that the map $\P A \to \P X$ is a trivial cofibration of $\mathbf{L}_{\mathbf{S}^{n} \subset \Di^{n+1}}\top$. The proof is complete using \cite[Proposition~1.5.4]{ref_model2}.
\epf

Theorem~\ref{fix} is sufficient to prove \cite[Theorem~V.3.5]{model2} (it is its only use) by using the q-cofibration $\mathbf{S}^n \subset \Di^{n+1}$ to force $\pi_n({\P} U_n)\longrightarrow \pi_n({\P} X)$ to become one-to-one, and the q-cofibration $\{0\} \subset \mathbf{S}^{n+1}$, which is a pushout of the preceding one, to force $\pi_{n+1}({\P} U_n,\gamma)\longrightarrow \pi_{n+1}({\P} X,\gamma)$ to become onto\footnote{A better terminology would be to use the term CW-cofibrant instead of strongly cofibrant; the proof of \cite[Theorem~V.3.5]{model2} mimicks the usual proof to build a CW-approximation.}.

Finally, the proof of \cite[Theorem~III.5.2]{model2} is not complete either. Let us start by recalling \cite[Proposition~III.5.1]{model2} (the compactness hypothesis on $U$ is removed; it is assumed in \cite{model2} but it is useless; only the compactness of the segment $[0,1]$ matters): 

\bp \label{rewording}
Let $X$ be a cellular object of the q-model category of multipointed $d$-spaces~\footnote{Such an object is called a globular complex in \cite{model2}; I introduced the notion of multipointed $d$-space only 4 years later.} Let $\P X$ be the quotient of $\P^{\mathcal{G}}X$ by the action of $\mathcal{G}$. Let $p_X:\P^{\mathcal{G}}X\to \P X$ be the canonical map. Let $U$ be an object of $\top$. Let $f,g:U\to \P^{\mathcal{G}}X$ be two continuous maps such that $p_X.f=p_X.g$. Then there exists a unique map $\phi:U\to \mathcal{G}$ such that $f(u)=g(u).\phi(u)$ for all $u\in U$. Moreover, $\phi$ is necessarily continuous.
\ep

Proposition~\ref{rewording} holds in $\top$ because we entirely work in this proof in the underlying space of $X$ which is a Hausdorff $\De$-generated space and because, in all cases, $\P X$ is the quotient of the underlying set of $\P^{\mathcal{G}}X$ by the action of $\mathcal{G}$ equipped with the final topology since the action of $\mathcal{G}$ is continuous.

\bth \label{existence-section} (replacement for \cite[Theorem~III.5.2]{model2}) Let $X$ be a cellular object of the q-model category of multipointed $d$-spaces. Then the canonical map $p_X:\P^{\mathcal{G}}X\to \P X$ has a section $i_X:\P X\to \P^{\mathcal{G}}X$. \eth

The proof works in $\top$. The only thing that matters is that the underlying set of $\P X$ is the quotient of the underlying set of $\P^{\mathcal{G}}X$ by the action of $\mathcal{G}$. And in all cases, the topology of $\P X$ is the final topology since the action of $\mathcal{G}$ is continuous. 

The proof uses the fix of \cite[Proposition~15.1]{model3} expounded in Theorem~\ref{calculation-pathspace-flow}. A subsequent paper will use the theory of Moore flows to give a more conceptual proof of this theorem independent from \cite{model2}.

\bpf
The idea is still to build $i_X:\P X \to \P^{\mathcal{G}}X$ by transfinite induction on the cellular decomposition of $X$.  Suppose that we have the pushout diagram of flows
\[
\xymatrix@C=4em@R=4em{
	\glob(\de Z) \fr{g}\fd{} & A \fd{f} \\
	\glob(Z) \fr{} & \cocartesien X }
\]
such that the map $\de Z\to Z$ is a generating q-cofibration. Suppose that we have proved the existence of $i_A:\P A \to \P^{\mathcal{G}}X$. We are going to build $i_X:\P X\to \P^{\mathcal{G}}X$ by building a cocone \[\xi:(\D^f\stackrel{\bullet} \longrightarrow \P^{\mathcal{G}}X).\] The proof will be complete thanks to Theorem~\ref{path-almost-accessible}.

Consider the topological space $T$ defined by the pushout diagram of topological spaces 
\[
\xymatrix@C=4em@R=4em{
	\de Z \fr{g}\fd{} & \P_{g(0),g(1)}A \fd{f} \\
	Z \fr{} & \cocartesien T. }
\]
There exists a continuous map $h:Z\to [0,1]$ such that $h^{-1}(0)=\de Z$. By the universal property of the pushout, it is extended to a continuous map $h:T \to [0,1]$ such that $h^{-1}(0)= \P_{g(0),g(1)}A$. 

Let us introduce an equivalence relation on the objects of $\mathcal{P}^{g(0),g(1)}(A^0)$ as follows: 
	\[\underline{m}\simeq \underline{n}\hbox{ if and only if }\mathbb{S}(\underline{m}) = \mathbb{S}(\underline{n}).\]
We obtain a partition of the set of objects of $\mathcal{P}^{g(0),g(1)}(A^0)$ by this equivalence relation. Once the map $\xi_{\underline{n}}:  \D^f(\underline{n}) \to \P^{\mathcal{G}}X$ is constructed for a given non simplifiable object $\underline{n}$, the definition of $\xi_{\underline{m}}:  \D^f(\underline{m}) \to \P^{\mathcal{G}}X$ is forced on any $\underline{m}$ such that $\underline{m} \simeq \underline{n}$.

We are going to build the maps $\xi_{\underline{n}}:  \D^f(\underline{n}) \to \P^{\mathcal{G}}X$ by proceeding by induction on the height of the non-simplifiable object $\underline{n}$. There is nothing to do for the height $0$: the non-simplifiable objects of height $0$ are all tuples $((\alpha,0,\beta))$ for $\alpha$ and $\beta$ running over $A^0$.

Let us expound the induction on a particular case: it is always the same method. Let $u=g(0)$ and $v=g(1)$. We suppose that the map $\xi_{\underline{p}}:\D^f(\underline{p}) \longrightarrow  \P^{\mathcal{G}}X$ is already constructed on all non simplifiable $\underline{p}$ of height at most $1$, and on all objects belonging to the same equivalence classes, and that the diagram 
\[
\xymatrix@C=4em@R=4em
{ \D^f(\underline{p}) \ar@{->}[d] \ar@{->}[r] & \ar@{->}[d] \D^f(\underline{q}) \\
\P^{\mathcal{G}}X \ar@{=}[r]	& \P^{\mathcal{G}}X
}
\]
is commutative for all maps $\underline{p}\to \underline{q}$ such that the vertical maps are already constructed. We want to construct, for example, the map 
\[
\xi_{((u,1,v),(v,0,u),(u,1,v))}:T\p \P_{v,u}A\p T = \D^f((u,1,v),(v,0,u),(u,1,v)) \longrightarrow \P^{\mathcal{G}}X.
\]
We consider the commutative diagram of topological spaces 
\[
\xymatrix@C=9em@R=4em
{
L_{((u,1,v),(v,0,u),(u,1,v))}\D^f \fd{} \fr{L_{((u,1,v),(v,0,u),(u,1,v))}\xi} &  \P^{\mathcal{G}}X \fd{p_X}\\
T \p \P_{v,u}A \p T  \fr{k} \ar@{-->}[ru]^-{} & \P X.
}
\]
The top arrow $L_{((u,1,v),(v,0,u),(u,1,v))}\xi$ is constructed by induction hypothesis. The bottom arrow $k$ is induced by the normalized composition law $*_N$ of the multipointed $d$-space $X$. By induction hypothesis, there is a map $\xi_{((u,1,v))}:T \to \P^{\mathcal{G}}X$ and a continuous map $\xi_{((v,0,u))}:\P_{v,u}A \to \P^{\mathcal{G}}X$. Consider the map $\mathcal{L}:T \p \P_{v,u}A \p T \to \P^{\mathcal{G}}X$ defined by~\footnote{The choice $1/3,1/3,1/3$ is not important, the proof works by choosing any $\ell_1,\ell_2,\ell_3>0$ such that $\ell_1+\ell_2+\ell_3=1$.} \[\mathcal{L}(z_1,\gamma,z_2) = \big[\big(\xi_{((u,1,v))}(\delta_{z_1})\big)\mu_{1/3}\big] * \big[\big(\xi_{((v,0,u))}(\gamma)\big)\mu_{1/3}\big] * \big[\big(\xi_{((u,1,v))}(\delta_{z_2})\big)\mu_{1/3}\big].\] By construction, we have $p_X\mathcal{L} = k$. However, there is no reason for the top triangle to be commutative as well. If $z_1\in \de Z$ or $z_2\in \de Z$, then the execution paths $\mathcal{L}(z_1,\gamma,z_2)$ and $(L_{((u,1,v),(v,0,u),(u,1,v))}\xi)(z_1,\gamma,z_2)$ have the same image $z_1*\gamma*z_2$ by $p_X:\P^{\mathcal{G}}X\to \P X$. By Proposition~\ref{rewording} applied to $U=L_{((u,1,v),(v,0,u),(u,1,v))}\D^f$, there exists a (unique) continuous map \[\phi:L_{((u,1,v),(v,0,u),(u,1,v))}\D^f \to \mathcal{G}\] such that \[\mathcal{L}(z_1,\gamma,z_2)\phi(z_1,\gamma,z_2) = \big(L_{((u,1,v),(v,0,u),(u,1,v))}\xi\big)(z_1,\gamma,z_2).\] Then we set 
\begin{multline*}
\xi_{((u,1,v),(v,0,u),(u,1,v))}(z_1,\gamma,z_2)\\:= \mathcal{L}(z_1,\gamma,z_2)\big((1-\min(h(z_1),h(z_2)))\phi(z_1,\gamma,z_2) + \min(h(z_1),h(z_2))\id_{[0,1]}\big).
\end{multline*}
The definition above gives a well-defined map because the barycenter \[(1-\min(h(z_1),h(z_2)))\phi(z_1,\gamma,z_2) + \min(h(z_1),h(z_2))\id_{[0,1]}\] belongs to $\mathcal{G}$. And we extend the definition of $\xi_{\underline{n}}$ to all objects $\underline{n}$ of $\mathcal{P}^{g(0),g(1)}(A^0)$ such that $\underline{n}\simeq ((u,1,v),(v,0,u),(u,1,v))$.

We have to verify that each diagram like 
\begin{equation}
\xymatrix@C=4em@R=4em
{ \D^f(\underline{p}) \ar@{->}[d]_-{\xi_{\underline{p}}} \ar@{->}[r] & \ar@{->}[d]^-{\xi_{\underline{q}}} \D^f(\underline{q}) \\
	\P^{\mathcal{G}}X \ar@{=}[r]	& \P^{\mathcal{G}}X
}
\label{square}
\end{equation}
where the vertical maps are constructed is commutative. 

\begin{figure}
	\[
	\xymatrix@C=4em@R=4em
	{ \D^f(\underline{p}) \ar@{->}[d] \ar@{->}[r] \ar@/_40pt/@{->}[dd]_-{\xi_{\underline{p}}}& \ar@{->}[d] \D^f(\underline{q}) \ar@/^40pt/@{->}[dd]^-{\xi_{\underline{q}}}\\
		\D^f(\mathbb{S}(\underline{p}))\ar@{=}[r]\ar@{->}[d]	& \D^f(\mathbb{S}(\underline{q}))\ar@{->}[d]\\
		\P^{\mathcal{G}}X \ar@{=}[r]	& \P^{\mathcal{G}}X.
	}
	\]
	\caption{Decomposition of Diagram~\ref{square}}
	\label{dec1}
\end{figure}

If the map $\underline{p}\to \underline{q}$ belongs to $\mathcal{P}^{g(0),g(1)}(A^0)_-$, then $\mathbb{S}(\underline{p}) = \mathbb{S}(\underline{q})$ and Diagram~\ref{square} can be decomposed as in Figure~\ref{dec1}. The top square is commutative because by Proposition~\ref{final-matching}, $\mathbb{S}(\underline{p}) = \mathbb{S}(\underline{q})$ is the terminal object of the equivalence class. The bottom square is trivially commutative. It is in fact something general: there is never nothing to verify if the top map $\underline{p}\to \underline{q}$ belongs to $\mathcal{P}^{g(0),g(1)}(A^0)_-$.

\begin{figure}
\[ \xymatrix@C=4em@R=2em
	{
		\D^f(\underline{p}) \ar@/_40pt/@{->}[ddd]_-{\xi_{\underline{p}}}\fd{-}\fR{+} \ar@{}[drr] |{\circled{A}} &&\D^f(\underline{q}) \fd{-}\ar@/^40pt/@{->}[ddd]^-{\xi_{\underline{q}}}\\
		\D^f(\underline{m}) \ar@{->}[rd]^-{+}\fd{-} \ar@{->}[rr]^-{+} \ar@{}[drr] |{\circled{B}} && \D^f(\mathbb{S}(\underline{q})) \fD{}\\
		\D^f(\mathbb{S}(\underline{m})) \fd{}\ar@{}[dr] |{\circled{D}} & L_{\underline{q}}\D^f  \ar@{->}[ru]^-{+}\ar@{->}[rd]\ar@{}[r] |-{\circled{C}}  &   \\
		\P^{\mathcal{G}}X \ar@{=}[rr] && \P^{\mathcal{G}}X.
	}\]
	\caption{Decomposition of Diagram~\ref{square}}
	\label{dec2}
\end{figure}

Suppose now that the map $\underline{p}\to \underline{q}$ belongs to $\mathcal{P}^{g(0),g(1)}(A^0)_+$. We have to verify the commutativity of the diagram with the newly defined maps. Diagram~\ref{square} can be decomposed as in Figure~\ref{dec2}. The existence of $\underline{m}$ and the commutativity of $\circled{A}$ comes from the Reedy structure of $\mathcal{P}^{g(0),g(1)}(A^0)$. The commutativity of the triangle $\circled{B}$ comes from the universal property of the colimit. The commutativity of the triangle $\circled{C}$ comes from the method for constructing $\D^f(\mathbb{S}(\underline{q})) \to \P^{\mathcal{G}}X$. Finally, the commutativity of the triangle $\circled{D}$ is the induction hypothesis when the map $\D^f(\mathbb{S}(\underline{q})) \to \P^{\mathcal{G}}X$ is constructed. 
\epf

\subsection*{About the category of topological spaces chosen in the past papers.} All past papers can be adapted to $\Htop_{\mathbf{\De}}$, including \cite{model3} which is even a little bit simpler because $\Htop_{\mathbf{\De}}$ is locally presentable: the smallness condition becomes trivial. It is not clear that the proofs of \cite{model3} are valid without a separation condition. Indeed, the closedness of some diagonal is used in \cite[Proposition~10.5]{model3}. Valid proofs of the existence of the q-model structure of flows with or without a separation condition can be found in \cite[Theorem~3.11]{leftdetflow} and \cite[Theorem~7.4]{QHMmodel}.

\appendix

\section{Basic properties of the category of all diagrams}
\label{alldg}

Let $\K$ be a bicomplete category. We gather some basic results about the category $\D\K$ of all small diagrams over all small categories defined as follows. 

An object is a functor $F:{I}\to \K$ from a small category ${I}$ to $\K$. A morphism from $F:{I}_1\to \K$ to $G:{I}_2\to \K$ is a pair $(f:{I}_1\to {I}_2,\mu:F \Rightarrow G.f)$ where $f$ is a functor and $\mu$ is a natural transformation. If $(g,\nu)$ is a map from $G:{I}_2\to \K$ to $H:\underline{K}\to \K$, then the composite $(g,\nu).(f,\mu)$ is defined by $(g.f,(\nu.f)\odot\mu)$. The identity of $F:{I}_1\to \K$ is the pair $(\id_{{I}_1},\id_F)$. If $(h,\xi):(H:{K}\to \K) \to (I:{L} \to \K)$ is another map of $\D\K$, then we have
\begin{align*}
\left((h,\xi).(g,\nu)\right).(f,\mu)&=  (h.g,\xi.g \odot\nu).(f,\mu) \\
& = (h.g.f,\xi.g.f\odot\nu.f\odot\mu)\\
& = (h,\xi). (g.f,\nu.f\odot \mu)\\
& = (h,\xi) .\left((g,\nu).(f,\mu)\right).
\end{align*}
Thus the composition law is associative and the category $\D\K$ is well-defined.

\bp The forgetful functor $p:\D\K\to \dcat$ is a bifibred category. The category $\D\K$ is bicomplete. The forgetful functor $p:\D\K\to \dcat$ is limit-preserving and colimit-preserving. \ep

I learnt the fibred category argument from \cite{271951} and from a remark after the question \cite{266597}.

\bpf
Consider the forgetful functor $p:\D\K\to \dcat$ taking a diagram $F:{I}\to \K$ to the small category $\underline{I}$. The functor category $\K^{{I}}$ is the fibre of $p$ over ${I}$. Let $f:{I} \to {J}$ be a functor between small categories. For a functor $G:{J}\to \K$, let $f^*(G)=G.f$. There is a canonical map $f^*(G)\to G$ in $\D\K$ defined by the pair \[(f,\id_{G.f}:G.f\Longrightarrow G.f).\] For a functor $F:{I}\to \K$ , let $f_*(F)=\Lan_fF$. Since $\K$ is bicomplete by hypothesis, there is an adjunction $f_*\dashv f^*$ between $\K^{{I}}$ and $\K^{{J}}$. Let $F:{I}\to \K$ and $G:{J}\to \K$ be two objects of $\D\K$. Let $\omega=(f,\mu):F\to G$ be a map of $\D\K$. A factorization of $\omega$ as a composite in $\D\K$ (with the left-hand map vertical) \[\xymatrix@C=5em
{
	F \fr{\omega^f=(\id_{\underline{I}},\overline{\omega})} & f^*(G) \fr{} & G
}\]
implies $\id_{G.f}.\overline{\omega}=\mu$. We obtain $\overline{\omega}=\mu$ as the unique possible choice. Let $g:{J}\to {K}$ be another map of $\dcat$.  One has for all functors $H:{K}\to\K$ the isomorphisms of functors \[f^*(g^*(H))\iso H.g.f\iso(g.f)^*(H).\] By \cite[Proposition~3.1]{QHMmodel}, the forgetful functor $p:\D\K\to \dcat$ is a bifibred category. Every fibre over a small category is a category of diagrams over a fixed small category: therefore all fibres of the bifibred category $p:\D\K\to\dcat$ are bicomplete. Moreover, the category of small categories is bicomplete as well. Using \cite[Proposition~3.3]{Roig} and \cite[Proposition~3.3~\textdegree]{Roig}, we deduce that $\D\K$ is bicomplete. The fact that the forgetful functor $p:\D\K\to \dcat$ is limit-preserving and colimit-preserving comes from \cite[Proposition~3.3]{Roig} and \cite[Proposition~3.3\textdegree]{Roig}.
\epf

\bp \label{colimdk} The colimit functor induces a well-defined functor $\liminj : \D\K\to \K$ which is a left adjoint. \ep

\bpf Consider the functor $\mathcal{I}:\K\to \D\K$ which takes an object $X$ of $\K$ to the constant diagram $\Delta_{\mathbf{1}}X$ over the terminal small category $\mathbf{1}$. Then we have the sequence of natural isomorphisms (where $F:{I}\to \K$ is an object of $\D\K$)
\begin{align*}
\D\K(F,\mathcal{I}(X)) &\iso \D\K(F,\Delta_{\mathbf{1}}X) & \hbox{by definition of $\mathcal{I}$}\\
& \iso \K^{{I}}(F,f^*(\Delta_{\mathbf{1}}X)) & \hbox{where $f:{I}\to\mathbf{1}$ is the canonical functor}\\
& \iso \K^{{I}}(F,\Delta_{{I}}X) & \hbox{by definition of $f$}\\
&\iso \K(\liminj F,X) &\hbox{by adjunction.}
\end{align*}
This sequence of natural isomorphisms implies that the mapping $\liminj:\D\K\mapsto\K$ yields a well-defined functor and that it is a left adjoint.
\epf

In particular we have for the sum the following result: 

\bp \label{plusdk} Let $(\C_i)_{i \in I}$ be a small family of small categories. Let $D_i:\C_i\to \K$ be a functor for all $i\in I$. Then the diagram $\bigsqcup_{i \in I} D_i$ (in $\D\K$) is isomorphic to the unique diagram from $\bigsqcup_{i \in I} \C_i$ to $\K$ such that the restriction to $\C_i$ is $D_i$. 
\ep

\bpf Obvious.
\epf

Finally, we also need the following result: 

\bp \label{cartesiandk} Suppose that the bicomplete category $\K$ is cartesian closed. Then the colimit functor $\liminj : \D\K\to \K$ commutes with binary products.
\ep

\bpf
Write $[-,-]$ for the internal hom of $\K$. Let $A:I\to \K$ and $B:J\to \K$ be two objects of $\D\K$. Then we have the sequence of isomorphisms 
\begin{align*}
\K((\liminj A)\p (\liminj B),C) &\iso \K(\liminj A,[\liminj B,C])\\
&\iso \limproj_{i\in I} \K(A_i,[\liminj B,C])\\
&\iso \limproj_{i\in I} \K(\liminj B,[A_i,C]) \\
&\iso \limproj_{i\in I} \limproj_{j\in J} \K(A_i\p B_j,C)\\
&\iso \limproj_{(i,j)\in I\p J} \K(A_i\p B_j,C)\\
&\iso \K(\liminj(A\p B),C),
\end{align*}
for all objects $C$ of $\K$, the first and the third isomorphisms because $\K$ is cartesian closed, the fifth one by the definition of an inverse limit in the category of sets, and the second, the fourth and the last one by the universal property of the limits. The proof is complete thanks to the Yoneda lemma. 
\epf

\section{\texorpdfstring{$\De$}{Lg}-Hausdorff \texorpdfstring{$\De$}{Lg}-generated spaces}
\label{HTop}

This appendix expounds the definition, the cartesian closedness, the calculation of some colimits, the local presentability and finally the model structures. We adapt in the sequel the proofs found in \cite{Ref_wH}, \cite{Ref_wH_Strick} and \cite{Ref_wH_Rezk}.

\subsection*{Definition} We refer to \cite[Chapter~VI]{topologicalcat} or \cite[Chapter~7]{Borceux2} for the notion of topological functor. The object of the category  $\top_{\mathbf{\De}}$ are called the \textit{$\De$-generated spaces}. Using space-filling curves, one sees immediately that this category contains all disks, all cubes, all spheres and all simplices. For a tutorial about these topological spaces, see for example \cite[Section~2]{mdtop}. The category  $\top_{\mathbf{\De}}$ is locally presentable by \cite[Corollary~3.7]{FR} and cartesian closed. The internal hom functor is denoted by $\ttop(-,-)$. We denote by $\omega:\mathcal{T\!O\!P}\to \set$ the underlying set functor where $\mathcal{T\!O\!P}$ is the category of general topological spaces. It is fibre-small and topological. The restriction functor $\omega:\top_{\mathbf{\De}}\subset \mathcal{T\!O\!P}\to \set$ is fibre-small and topological as well. The category  $\top_{\mathbf{\De}}$ is a full coreflective subcategory of the category $\mathcal{T\!O\!P}$ of general topological spaces. Let $k_\De:\mathcal{T\!O\!P}\to\top_{\mathbf{\De}}$ be the $\De$-kelleyfication functor, i.e. the right adjoint. The category  $\top_{\mathbf{\De}}$ is finally closed in $\mathcal{T\!O\!P}$, which means that the final topology and the $\omega$-final structure coincides. On the contrary, the $\omega$-initial structure in  $\top_{\mathbf{\De}}$ is obtained by taking the $\De$-kelleyfication of the initial topology in $\mathcal{T\!O\!P}$. If $A$ is a subset of a space $X$ of  $\top_{\mathbf{\De}}$, the initial structure in  $\top_{\mathbf{\De}}$ of the inclusion $A\subset \omega X$ is the $\De$-kelleyfication of the relative topology with respect to the inclusion.

\bd A general topological space $X$ is {\rm $\De$-Hausdorff} if for every continuous map $f:[0,1] \to X$, the subset $f([0,1])$ is closed in $X$. \ed

In particular, every point of a $\De$-Hausdorff general topological space is closed (i.e. every $\De$-Hausdorff general topological space is a $T_1$-space) and every finite subset of a $\De$-Hausdorff general topological space equipped with the relative topology is discrete.

\bp \label{HwH} A Hausdorff general topological space is $\De$-Hausdorff. \ep

\bpf Let $X$ be a Hausdorff general topological space. Let $f:[0,1] \to X$ be a continuous map. Then $f([0,1])$ is quasi-compact, i.e. it satisfies the finite open covering property and it is Hausdorff since $X$ is Hausdorff. Thus $f([0,1])$ is compact and therefore closed since $X$ is Hausdorff.
\epf

\bp \label{carac} A $\De$-generated space $X$ is $\De$-Hausdorff if and only if the diagonal $\Delta_X=\{(x,x)\mid x\in X\}$ is a closed subset of $X\p X$ where the product is taken in $\top_{\mathbf{\De}}$. \ep

\bpf Suppose that $\Delta_X$ is a closed subset of $X\p X$. Let $f:[0,1]\to X$ be a continuous map. Let $v:[0,1]\to X$ be another continuous map. Let $M=\{(a,b)\in [0,1]\p [0,1]\mid f(a)=v(b)\}\subset [0,1]\p [0,1]$. Since $M=(f\p v)^{-1}(\Delta_X)$, it is closed in $[0,1]\p [0,1]$, thus compact. We deduce that the projection $\pi_{[0,1]}(M)$ is compact and thus closed in $[0,1]$. However $\pi_{[0,1]}(M)=v^{-1}(f([0,1]))$. We deduce that $f([0,1])$ is closed in $X$. Conversely, assume that $X$ is $\De$-Hausdorff. Then every one-point subset $\{x\}\subset X$ is closed in $X$ as the image of a constant map from $[0,1]$ to $X$. Consider a map $u=(v,w):[0,1] \to X\p X$. It is enough to show that $u^{-1}(\Delta_X)=\{a\in [0,1] \mid v(a)=w(a)\}$ is closed in $[0,1]$. Suppose that $a\notin u^{-1}(\Delta_X)$, ie. $v(a)\neq w(a)$. Since $\{w(a)\}$ is closed in $X$, the set $U=\{b\mid v(b)\neq w(a)\}$ is an open subset of $[0,1]$ containing $a$. There is an open interval $]a-\epsilon,a+\epsilon[$ of $[0,1]$ containing $a$ such that $[a-\epsilon,a+\epsilon]\subset U$, or equivalently $w(a) \notin v([a-\epsilon,a+\epsilon])$. This implies that $a \in W=w^{-1}(X\backslash v([a-\epsilon,a+\epsilon]))$. Since $X$ is $\Delta$-Hausdorff, the set $v([a-\epsilon,a+\epsilon])$ is a closed subset of $X$. Thus $W$ is an open subset of $[0,1]$. Let $b\in ]a-\epsilon,a+\epsilon[ \cap W$. Then $v(b)\in v([a-\epsilon,a+\epsilon])$ and $w(b)\notin v([a-\epsilon,a+\epsilon])$. Thus $b\notin u^{-1}(\Delta_X)$. It means that $]a-\epsilon,a+\epsilon[ \cap W$ is an open subset of $[0,1]$ containing $a$ included in the complement of $u^{-1}(\Delta_X)$. Thus $u^{-1}(\Delta_X)$ is closed and so is $\Delta_X$. 
\epf

\subsection*{Cartesian closedness}

\bp \label{wprod} The product in $\mathcal{T\!O\!P}$ of an arbitrary family of $\De$-Hausdorff general topological spaces $X_i$ with $i$ running over a set of indices $I$ is $\De$-Hausdorff. \ep

\bpf
Let $X=\prod_{i\in I} X_i$ with projection maps $\pr_i:X\to X_i$. Let $f:[0,1] \to X$ be a continuous map. Then for all $i\in I$, the set $\pr_i(f([0,1]))$ is closed in $X_i$. Thus, $X$ being equipped with the initial topology making the projection maps continuous, $f([0,1])$ is closed in $X$. We deduce that $X$ is $\De$-Hausdorff. 
\epf

\bp \label{wincl} Let $i:A\to B$ be a one-to-one continuous map between $\De$-generated spaces. If $B$ is $\De$-Hausdorff, then $A$ is $\De$-Hausdorff as well. In particular, if $A$ is a subset of $B$ equipped with the relative topology, then if $B$ is $\De$-Hausdorff, then $A$ is $\De$-Hausdorff as well.
\ep

\bpf
Let $p:[0,1]\to A$ be a continuous map. Then $(ip)([0,1])$ is closed in $B$ because $B$ is $\De$-Hausdorff by hypothesis. Therefore the inverse image $i^{-1}((ip)([0,1]))$ is closed in $A$. Since $i$ is one-to-one, we deduce that $i^{-1}((ip)([0,1]))=p([0,1])$. We obtain that $A$ is $\De$-Hausdorff. 
\epf

\bp \label{internalhomhausdorff} Let $X$ and $Y$ be two $\De$-generated spaces with $Y$ $\De$-Hausdorff. Then $\ttop(X,Y)$ is $\De$-Hausdorff. \ep

\bpf
There is one-to-one continuous map $g:\ttop(X,Y) \to \prod_{x\in X} Y$ induced by the evaluation maps, the product being taken in  $\top_{\mathbf{\De}}$. By Proposition~\ref{wprod} and since the $\De$-kelleyfication adds closed subsets, the space $\prod_{x\in X} Y$ is $\De$-Hausdorff. The proof is complete thanks to Proposition~\ref{wincl}. 
\epf

\bp \label{refl} The category  $\Htop_{\mathbf{\De}}$ of $\De$-Hausdorff $\De$-generated spaces is a full reflective subcategory of  $\top_{\mathbf{\De}}$. \ep

\bpf[Sketch of proof]
The $\De$-Hausdorffization functor $w_\Delta:\top_{\mathbf{\De}}\to \Htop_{\mathbf{\De}}$ looks like the  weakly Hausdorffization functor. Starting from a $\De$-generated space $X$, take two points $x$ and $y$ such that $(x,y)$ belongs to the closure of $\De_X$ and consider the quotient $X\to X/(x=y)$ with the final topology. Iterate the process transfinitely. It will stop eventually for a cardinality reason. The proof is then formally the same as in the case of weakly Hausdorff $k$-spaces. 
\epf

\begin{cor} \label{HCartClosed} The category  $\Htop_{\mathbf{\De}}$ of $\De$-Hausdorff $\De$-generated spaces is cartesian closed. \end{cor}

\bpf
The product in  $\top_{\mathbf{\De}}$ of two $\De$-Hausdorff $\De$-generated spaces is $\De$-Hausdorff by Proposition~\ref{refl}. The proof is complete thanks to Proposition~\ref{internalhomhausdorff}. 
\epf

\subsection*{Calculation of some colimits}

\bd A {\rm quotient map} is a continuous map $f:X\to Y$ of $\top_{\mathbf{\De}}$ which is onto and such that $Y$ is equipped with the final topology. The space $Y$ is called a {\rm final quotient} of $X$.
\ed

If $f:X\to Y$ is a surjective continuous map of $\top_{\mathbf{\De}}$ which is either open or closed, it is easy to see that it is a quotient map.

\bp \label{quotient-and-Hausdorff} Let $f:X\to Y$ be a quotient map of $\De$-generated spaces. Then $Y$ is $\De$-Hausdorff if and only if the set $\{(x,y)\in X\p X\mid f(x)=f(y)\}$ is a closed subset of $X\p X$
\ep

\bpf One has $\{(x,y)\in X\p X\mid f(x)=f(y)\} = (f\p f)^{-1}(\Delta_Y)$. If $Y$ is $\De$-Hausdorff, then $\Delta_Y$ is closed in $Y\p Y$ by Proposition~\ref{carac} and therefore $\{(x,y)\in X\p X\mid f(x)=f(y)\}$ is closed in $X\p X$. Conversely, if $\{(x,y)\in X\p X\mid f(x)=f(y)\}=(f\p f)^{-1}(\Delta_Y)$ is closed then $\Delta_Y$ is closed because $Y\p Y$ is equipped with the final topology. Thus by Proposition~\ref{carac}, $Y$ is $\De$-Hausdorff. 
\epf

\begin{cor} \label{closed-eq}
	Let $E$ be an equivalence relation on a $\De$-generated space $X$. Then $X/E$ equipped with the final topology is $\De$-Hausdorff if and only if $E$ is a closed subset of $X\p X$.  
\end{cor}

\begin{cor} \label{ClosedSubspace}
	Let $A$ be a closed subset of a $\De$-Hausdorff $\De$-generated space. Then the quotient $X/A$ equipped with the final topology is $\De$-Hausdorff. 
\end{cor}

Unlike the case of $k$-spaces, a closed subset of a $\De$-generated space equipped with the relative topology is not necessarily $\De$-generated. For example, the Cantor set $\mathbb{K}\subset [0,1]$ is closed and totally disconnected. Therefore, every continuous map from $[0,1]$ to $\mathbb{K}$ is constant, which implies that its $\De$-kelleyfication is a discrete set. Thus $\mathbb{K}$ is not $\De$-generated whereas the segment $[0,1]$ is $\De$-generated.

\bd \label{def_closed_incl_delta} A {\rm closed inclusion} of $\De$-generated spaces $i:A\to B$ is a one-to-one continuous map of  $\top_{\mathbf{\De}}$ such that $i(A)$ is a closed subset of $B$ and such that $i$ induces a homeomorphism between $A$ and $i(A)$ equipped with the relative topology (which implies that $i(A)$ is $\De$-generated). 
\ed

\bp \label{L1} Consider the commutative diagram of  $\top_{\mathbf{\De}}$: 
\[
\xymatrix@C=4em@R=4em
{
	X \fr{p} \ar@{->}[d]_-{f} & Z \ar@{->}[d]^-{g} \\
	Y \fr{q} & W
}
\]
If $f$ is a closed inclusion, $g$ is one-to-one, $p$ is onto and either $q$ is closed (i.e. the image of a closed subset is a closed subset) or $q$ is a quotient map, with $q^{-1}(g(Z))\subset f(X)$, then $g$ is a closed inclusion.
\ep

\bpf
Let $F$ be a closed subset of $Z$. Then we have $q^{-1}(g(F))\subset q^{-1}(g(Z))\subset f(X)$. We have $f^{-1}(q^{-1}(g(F))) = p^{-1}(F)$: thus $q^{-1}(g(F))$ is a closed subset of $f(X)$ equipped with the relative topology because $f:X\to Y$ is a closed inclusion by hypothesis. Thus $q^{-1}(g(F))$ is a closed subset of $Y$ and since $q$ is a quotient map or a closed map, $g(F)$ is a closed subset of $W$. In particular, $g(Z)$ is a closed subset of $W$ and there is the homeomorphism $Z\iso g(Z)$ (which implies that $g(Z)$ equipped with the relative topology is $\De$-generated). 
\epf

\bp \label{push-closed-incl} Consider the pushout diagram of  $\top_{\mathbf{\De}}$
\[
\xymatrix@C=4em@R=4em
{
	A \fr{h} \fd{f}& X \ar@{->}[d]^-{g} \\
	B \fr{k}& \cocartesien Y
}
\]
If the map $f:A\to B$ is a closed inclusion of  $\top_{\mathbf{\De}}$, then the map $g:X\to Y$ is a closed inclusion of $\top_{\mathbf{\De}}$. If moreover $B$ and $X$ are $\De$-Hausdorff, then $Y$ is $\De$-Hausdorff. 
\ep

Note that if $B$ is $\De$-Hausdorff, then $A$ is $\De$-Hausdorff as well by Proposition~\ref{wincl}.

\bpf
The injectivity of $g$ comes from the fact that  $\top_{\mathbf{\De}}$ is topological over the category of sets. Consider the commutative diagram of  $\top_{\mathbf{\De}}$
\[
\xymatrix@C=4em@R=4em
{
	X\sqcup A \fr{(\id,h)} \fd{\id \sqcup f}& X \ar@{->}[d]^-{g} \\
	X\sqcup B \fr{(g,k)}&  Y
}
\]
Then $(g,k)^{-1}(g(X)) \subset X\sqcup k^{-1}(g(X)) \subset X\sqcup f(A)$, the last inclusion because the diagram of the statement of the proposition is a pushout. However, $(g,k)$ is a quotient map. So by Proposition~\ref{L1}, $g:X\to Y$ is a closed inclusion. 

The diagram above is a pushout diagram. Moreover, the map $\id \sqcup f$ is a closed inclusion and the map $(\id,h)$ is a quotient map. Therefore, we can suppose in the original diagram without lack of generality that $k$ is a quotient map. We have $(k\p k)^{-1}(\De_Y) = \De_\B \cup (f\p f).(h\p h)^{-1}(\De_X)$ since $f$ is one-to-one. Assume that $B$ and $X$ are $\De$-Hausdorff. Then $(h\p h)^{-1}(\De_X)$ is closed in $A\p A$ and therefore $(f\p f).(h\p h)^{-1}(\De_X)$ is closed in $B\p B$. And $\De_\B$ is closed in $B\p B$. Thus $(k\p k)^{-1}(\De_Y)$ is closed in $B\p B$. Using Proposition~\ref{carac}, we deduce that $Y$ is $\De$-Hausdorff. 
\epf

\bp \label{alreadyH} Let $\lambda$ be a limit ordinal. Consider a transfinite tower $X:\lambda\to \Htop_{\mathbf{\De}}$ of $\De$-Hausdorff $\De$-generated spaces such that for every ordinal $\alpha<\lambda$, the map $X_\alpha\to X_{\alpha+1}$ is one-to-one. Then the colimit $\liminj X$ calculated in  $\top_{\mathbf{\De}}$ is $\De$-Hausdorff.
\ep

Note that it is not assumed that the maps $X_\alpha\to X_{\alpha+1}$ are closed inclusions. 

\bpf
Let $U=\bigsqcup_{\alpha<\lambda} X_\alpha$. Then there exists a quotient map $q:U\to \liminj X$. Let \[E=\{(x,y)\in U\p U\mid q(x)=q(y)\}.\] Let $E = \bigsqcup_{\alpha,\beta} E_{\alpha,\beta}$ where $E_{\alpha,\beta}= E \cap (X_\alpha \p X_\beta)$ with $\alpha,\beta<\lambda$. Let $\lambda>\delta\geq \max(\alpha,\beta)$. Let $f_{\alpha,\delta}:X_\alpha \to X_{\delta}$ and $f_{\beta,\delta}:X_\beta \to X_{\delta}$. Since each map of the tower is one-to-one, the map $X_\delta \to \liminj X$ is one-to-one, and we have therefore $E_{\alpha,\beta} = (f_{\alpha,\delta}\p f_{\beta,\delta})^{-1}(\De_{X_\delta})$ for any $\lambda>\delta\geq \max(\alpha,\beta)$. We then deduce from the fact that $X_\delta$ is $\De$-Hausdorff and using Proposition~\ref{carac} that $E_{\alpha,\beta}$ is a closed subset of $U\p U$ for all $\alpha,\beta<\lambda$. It follows that $\liminj X$ is $\De$-Hausdorff using Corollary~\ref{closed-eq}. 
\epf

\subsection*{Local presentability}

\bp \label{carH2} A $\De$-generated space is $\De$-Hausdorff if and only if it has unique sequential limits. \ep

\bpf
Every $\De$-generated space is a sequential space because the segment $[0,1]$ is sequential. The proof is complete thanks to Proposition~\ref{carac}. 
\epf

\bp \label{HLocPre} The category  $\Htop_{\mathbf{\De}}$ is locally presentable. \ep

\bpf
By Proposition~\ref{carH2}, to encode $\De$-Hausdorff $\De$-generated spaces, it suffices, using \cite[Theorem~5.3]{MR629337} and \cite[Theorem~3.6]{FR}, to start from a small relational universal strict Horn theory $\mathcal{T}$ without equality encoding  $\top_{\mathbf{\De}}$ and to encode the fact that the limit of a sequence $x$ is unique. Consider the set bijection \[f:\bigg\{\frac{1}{n}\mid n\geq 1\bigg\}\cup \{0\} \stackrel{\iso}\longrightarrow \mathbb{N}\cup \{\infty\}\] 
such that $f(1/n)=n-1$ for all $n\geq 1$ and such that $f(0)=\infty$. Put on the set $\{1/n\mid n\geq 1\}\cup \{0\}$ the relative topology induced by the one of $[0,1]$. We obtain using the set bijection $f$ a topology on $\mathbb{N}\cup \{\infty\}$. Observe that a sequence $x$ converges to $y$ in $X$ if and only if the sequence $xy$ yields a continuous map from $\mathbb{N}\cup \{\infty\}$ to $X$. Consider a continuous map $g:[0,1]\to \mathbb{N}\cup \{\infty\}$. If $g$ is not the constant map $\infty$, then there exists $n\in \mathbb{N}$ such that $g^{-1}(n)$ is nonempty. Since $g^{-1}(n)$ is both open and closed in $[0,1]$ which is connected, we deduce the equality $g^{-1}(n)=[0,1]$. Thus the constant maps are the only continuous maps from $[0,1]$ to $\mathbb{N}\cup \{\infty\}$. Therefore the $\De$-kelleyfication of the topological space $\mathbb{N}\cup \{\infty\}$, which is equipped with the final topology with respect to all constant maps, is the discrete space with the same underlying set. It means that the topological space $\mathbb{N}\cup \{\infty\}$ is not $\De$-generated. We have therefore to consider the final closure $\mathrm{TOP}$ in $\mathcal{T\!O\!P}$ of $[0,1]$ and $\mathbb{N}\cup \{\infty\}$. By \cite[Theorem~5.3]{MR629337} and \cite[Theorem~3.6]{FR}, there exists another small relational universal strict Horn theory $\mathcal{T}'$ without equality encoding $\mathrm{TOP}$. Let $\{R_j\mid j\in J\}$ be the set of relational symbols of $\mathcal{T}'$. That the limit of a sequence is unique is then formalized by the axiom 
\[
(U)\ (\forall x,y_1,y_2) \big(\bigwedge_{(R_j,a)} R_j((xy_1).a) \wedge R_j((xy_2).a)\big) \longrightarrow y_1 = y_2
\]
where the conjonction is taken on all pairs $(R_j,a)$ such that $\mathcal{T}'$ satisfies $R_j$ for a sequence $a$ of elements of $\mathbb{N}\cup \{\infty\}$ ($(xy).a$ meaning the composition of $xy$ and $a$). The theory $\mathcal{T}\cup \mathcal{T}'$ with the axiom $(U)$ is a small relational universal strict Horn theory with equality which provides a model of  $\Htop_{\mathbf{\De}}$. By \cite[Theorem~5.30]{TheBook}, the category  $\Htop_{\mathbf{\De}}$ of $\De$-Hausdorff $\De$-generated spaces is therefore locally presentable. 
\epf

\begin{cor}
	The $\De$-Hausdorffization functor $w_\De:\top_{\mathbf{\De}}\to \Htop_{\mathbf{\De}}$ is accessible.
\end{cor}

\subsection*{The model structures}

Before introducing three model structures on  $\Htop_{\mathbf{\De}}$, let us recall Isaev's results adapted to our situation: 

\bth \cite[Theorem~4.3, Proposition~4.4, Proposition~4.5 and Corollary~4.6]{Isaev} \label{Isa}
Let $\K$ be a locally presentable category. Let $I$ be a set of maps of $\K$ such that the domains of the maps of $I$ are $I$-cofibrant (i.e. belong to $\cof(I)$). Suppose that for every map $i:U\to V \in I$, the relative codiagonal map $V\sqcup_U V \to V$ factors as a composite $V\sqcup_U V \stackrel{\gamma_0\sqcup \gamma_1}\to C_U(V)\to V$ such that the left-hand map belongs to $\cof(I)$. Let $J_I = \{\gamma_0:V\to C_U(V)\mid U\to V \in I\}$. Suppose that there exists a path functor $\cocyl:\K \to \K$, i.e. an endofunctor of $\K$ equipped with two natural transformations $\tau:\id\Rightarrow\cocyl$ and $\pi:\cocyl\Rightarrow \id\p \id$ such that the composite $\pi.\tau$ is the diagonal. Moreover we suppose that the path functor satisfies the following hypotheses: 
\begin{enumerate}
	\item With $\pi=(\pi_0,\pi_1)$, $\pi_0:\cocyl(X)\to X$ and $\pi_1:\cocyl(X)\to X$ have the RLP with respect to $I$.
	\item The map $\pi:\cocyl(X)\to X\p X$ has the RLP with respect to the maps of $J_I$.
\end{enumerate}
Then there exists a unique model category structure on $\K$ such that the set of generating cofibrations is $I$ and such that the set of generating trivial cofibrations is $J_I$. Moreover, all objects are fibrant.
\eth

Thanks to Corollary~\ref{HCartClosed} and Proposition~\ref{HLocPre}, we can use Theorem~\ref{Isa} to obtain on $\Htop_{\mathbf{\De}}$ the \textit{q-model structure} $(\C_q,\W_q,\F_q)$ by using the path space functor \[\cocyl:X\mapsto \ttop([0,1],X).\] The cofibrations, called \textit{q-cofibrations}, are the retracts of the transfinite compositions of the inclusions $\mathbf{S}^{n-1}\subset \Di^n$ for $n\geq 0$, the weak equivalences are the weak homotopy equivalences and the fibrations, called \textit{q-fibrations} are the maps satisfying the RLP with respect to the inclusions $\Di^{n}\subset \Di^{n+1}$ for $n\geq 0$, or equivalently with respect to the inclusions $\Di^{n}\p \{0\}\subset \Di^{n}\p [0,1]$ for $n\geq 0$; this model structure is combinatorial. It is Quillen equivalent to the q-model structure of  $\top_{\mathbf{\De}}$ because the canonical map $X\to w_\De (X)$ is an isomorphism for all q-cofibrant objects of  $\Htop_{\mathbf{\De}}$ by Proposition~\ref{HwH}.

Using Corollary~\ref{HCartClosed}, we can also put on  $\Htop_{\mathbf{\De}}$ a structure of topologically bicomplete category in the sense of \cite{Barthel-Riel}. Since the category  $\Htop_{\mathbf{\De}}$ is locally presentable by Proposition~\ref{HLocPre}, it satisfies the monomorphism hypothesis. Using \cite[Corollary~5.23]{Barthel-Riel}, we then obtain the \textit{h-model structure} $(\C_h,\W_h,\F_h)$: the fibrations, called the \textit{h-fibrations}, are the maps satisfying the RLP with respect to the inclusions $X\p \{0\}\subset X\p [0,1]$ for all topological spaces $X$, and the weak equivalences are the homotopy equivalences;  we have $\C_q\subset \C_h$ because $\W_h\subset \W_q$ and $\F_h\subset \F_q$. The cofibrations are called the \textit{h-cofibrations}. The h-cofibrations are the closed inclusions satisfying the LLP with respect to the maps of the form $\ttop([0,1],Y)\to Y$ \cite[Proposition~1(b)]{vstrom3} (a part of the argument is reproduced in Proposition~\ref{h-cof-relative-T1}). It is well-known that, in $\Htop_\mathbf{K}$, a map is a h-cofibration if and only it satisfies the LLP with respect to the map of the form $\ttop([0,1],Y)\to Y$. It is true as well in $\Htop_{\mathbf{\De}}$. If a map $i:A\to X$ satisfies the LLP with respect to the maps of the form $\ttop([0,1],Y)\to Y$, then the map $i\square i_0:Mi\to X\p [0,1]$ (with $i_0:\{0\}\subset [0,1]$) has a retract $r$ and therefore is a \textit{closed} inclusion because $Mi\iso \{y\in X\p [0,1]\mid ir(y)=y\}\subset X\p [0,1]$ is a closed subset since $X\p [0,1]$ is $\De$-Hausdorff. Then the classical argument applies \cite[Proposition~8.2]{Ref_wH} \cite[Lemma~1.6.2(ii)]{ParamHomTtheory}: consider the commutative diagram 
\[
\xymatrix@C=4em@R=4em
{
A \fd{i} \fr{a\mapsto (a,1)} & Mi \ar@{->}[d]^-{i\square i_0} \\
X   \fr{x\mapsto (x,1)} & X\p [0,1] .
}
\] 
Then $i$ is a closed inclusion because the three other ones are closed inclusions.

The \textit{m-model structure} $(\C_m,\W_m,\F_m)=(\C_m,\W_q,\F_h)$ is characterized as follow: the fibrations are the \textit{h-fibrations}, and the weak equivalences are the weak homotopy equivalences; we have $\C_q\subset \C_m$ because $\W_m\cap \F_m=\W_q\cap \F_h \subset \W_q\cap \F_q$. Its existence is a consequence of \cite[Theorem~2.1]{mixed-cole}. By \cite[Corollary~3.7]{mixed-cole}, a topological space is m-cofibrant if and only if it is homotopy equivalent to a q-cofibrant space. It is the mixed model structure in the sense of \cite{mixed-cole} of the two preceding model structures. We have $\C_m\subset \C_h$ by \cite[Proposition~3.6]{mixed-cole}.


\newcommand{\etalchar}[1]{$^{#1}$}

\end{document}